\documentclass[11pt]{amsart}

\usepackage{amssymb}
\usepackage{float}

\usepackage{fullpage}

\usepackage{graphicx}
\usepackage{comment}
\usepackage{longtable}
\usepackage{enumerate}
\usepackage{enumitem}
\usepackage{array}
\usepackage{calc}
\usepackage{multicol}
\usepackage{color}

\newtheorem{theorem}{Theorem}[section]
\newtheorem{lemma}[theorem]{Lemma}
\newtheorem{proposition}[theorem]{Proposition}

\theoremstyle{definition}
\newtheorem{remark}[theorem]{Remark}

\newcommand\bthree[1]{{#1}}
\newcommand\ctwo[1]{{#1}}
\newcommand\cthree[1]{{#1}}
\newcommand\dfour[1]{{#1}}
\newcommand\dfive[1]{{#1}}
\newcommand\dsix[1]{{#1}}
\newcommand\dseven[1]{{#1}}

\newcommand{\myPhi}{\ensuremath{\mathrm{\Phi}}}
\newcommand{\myPi}{\ensuremath{\mathrm{\Pi}}}
\newcommand{\myTheta}{\ensuremath{\mathrm{\Theta}}}

\newcommand{\deq}{\overset{\textup{def}}{=}}
\DeclareMathOperator{\Stab}{Stab}

\DeclareMathOperator{\Aut}{Aut}
\DeclareMathOperator{\SL}{SL}
\DeclareMathOperator{\GL}{GL}
\DeclareMathOperator{\SO}{SO}
\DeclareMathOperator{\GO}{GO}
\DeclareMathOperator{\Sp}{Sp}

\setlist[enumerate]{topsep=0pt}

\title{Reducible subgroups of exceptional algebraic groups}

\author{Alastair J.\ Litterick}
\address{A.J. Litterick, Fakult\"{a}t f\"{u}r Mathematik, Ruhr-Universit\"{a}t Bochum, Universit\"{a}tsstra{\ss}e 150, D-44780 Bochum, Germany, and Fakult\"{a}t f\"{u}r Mathematik, Universit\"{a}t Bielefeld, Postfach 100131, D-33501 Bielefeld, Germany.}
\email{ajlitterick@gmail.com}
\thanks{The first author is supported by the Alexander von Humboldt Foundation. Both authors thank the anonymous referee for their careful reading and numerous suggestions that have improved the paper.}
\author{Adam R.\ Thomas}
\address{A.R. Thomas, School of Mathematics, University of Bristol, Bristol, BS8 1TW, UK, and The Heilbronn Institute for Mathematical Research, Bristol, UK.}
\email{adamthomas22@gmail.com}


\begin{document}

\begin{abstract}
Let $G$ be a simple algebraic group over an algebraically closed field. A closed subgroup $H$ of $G$ is called $G$-completely reducible ($G$-cr) if, whenever $H$ is contained in a parabolic subgroup $P$ of $G$, it is contained in a Levi factor of $P$. In this paper we complete the classification of connected $G$-cr subgroups when $G$ has exceptional type, by determining the $L_{0}$-irreducible connected reductive subgroups for each simple classical factor $L_{0}$ of a Levi subgroup of $G$. As an illustration, we determine all reducible, $G$-cr semisimple subgroups when $G$ has type $F_4$ and various properties thereof. This work complements results of Lawther, Liebeck, Seitz and Testerman, and is vital in classifying non-$G$-cr reductive subgroups, a project being undertaken by the authors elsewhere.
\end{abstract}

\maketitle

\section{Introduction}

The subgroup structure of linear algebraic groups is a theory with applications throughout algebra and beyond. When working over $\mathbb{R}$ or $\mathbb{C}$, close connections with compact Lie groups give direct applications in geometry and physics. Over local fields and global fields, important applications are found in number theory, for instance in the Langlands program. Over finite fields and their algebraic closures, the subgroup structure has direct applications to the study of finite simple groups, for instance the subgroup structure of exceptional simple algebraic groups has been heavily used in understanding maximal subgroups of the corresponding finite groups of Lie type \cite{Litterick2018,Craven2017,MR1132853,MR2044850}.

Reductive subgroups of reductive groups are an important special instance of this problem. On the one hand, a result of Borel and Tits \cite{MR0294349} states that for a reductive group $G$, a connected subgroup of $G$ is either itself reductive, or contained `canonically' in a proper parabolic subgroup of $G$. Thus the problem reduces to understanding parabolic subgroups and reductive subgroups, for instance a maximal connected subgroup of $G$ is either parabolic or reductive. On the other hand, understanding reductive subgroups is equivalent to understanding affine homogeneous spaces for $G$, since for a closed subgroup $H \le G$, the coset space $G/H$ is an affine variety precisely when $H^{\circ}$ is reductive \cite{Richardson1977}.

Let $G$ be a connected reductive algebraic group over an algebraically closed field. Following \cite{Ser3}, a closed subgroup $H$ of $G$ is said to be \emph{$G$-completely reducible} ($G$-cr) if whenever $H$ is contained in a parabolic subgroup $P$ of $G$, it is contained in a Levi factor of $P$, and \emph{non-$G$-cr} otherwise. Similarly, $H$ is called \emph{$G$-irreducible} or just \emph{irreducible} if $H$ is not contained in any proper parabolic subgroup of $G$, and \emph{reducible} otherwise. The $G$-cr subgroups of $G$ are precisely the $L$-irreducible subgroups of $L$, as $L$ ranges over the Levi subgroups of $G$ \cite[Corollary 3.5]{MR2178661}, allowing $L = G$ to account for $G$-irreducible subgroups.

In the present work we classify the $L_{0}$-irreducible connected subgroups as $L_{0}$ varies over those simple algebraic groups of classical type which occur as factors of Levi subgroups of a simple algebraic group $G$ of exceptional type. Together with work of the second author and many others \cite{MR2707891,MR1466951,MR1329942,MR2044850,MR3075783,MR2604850,Tho1,Tho2,Tho3}, the results here effectively complete the classification of $G$-cr connected subgroups for simple $G$ of exceptional type. In particular, the cited papers of the second author determine all $G$-irreducible connected subgroups for $G$ of exceptional type, and so it remains for us to describe $L$-irreducible connected subgroups for $L$ involving factors of classical type.

The information on completely reducible subgroups derived here is also vital in the complementary problem of understanding non-$G$-cr subgroups. A fruitful approach to this problem, taken for instance in \cite{MR3075783,LitTho,LitTho2}, is as follows. For each parabolic subgroup $P$ of $G$, fix a Levi factor $L$ of $P$ and take the (known) list of $L$-irreducible subgroups of $L$. For each such subgroup $X$, consider the set of complements to the unipotent radical $R_{u}(P)$ in the semidirect product $R_{u}(P)X$. As $P$, $L$ and $X$ vary, these complements give rise to all non-$G$-cr subgroups of $G$, modulo the subtle but tractable problem of understanding the difference between `$G$-conjugacy' and `$P$-conjugacy'.

By a `Levi subgroup of $G$' we mean a Levi factor of a parabolic subgroup of $G$. To state our main theorem, recall that for each $G$-cr subgroup $X \le G$, there exists a subgroup $L$ of $G$ which is minimal among Levi subgroups of $G$ containing $X$, and the conjugacy class of $L$ in $G$ is uniquely determined by the class of $X$ (cf.\ Proposition \ref{prop:glconj}).

\begin{theorem} \label{thm:main}
Let $G$ be a simple algebraic group of exceptional type over an algebraically closed field, let $X$ be a $G$-completely reducible connected subgroup of $G$, and let $L$ be minimal among Levi subgroups of $G$ containing $X$. Then the image of $X$ under projection to each classical simple factor of $L$ is given by an embedding in Section \ref{sec:embeddings}.
\end{theorem}

\begin{remark}
Irreducible subgroups of a simple algebraic group form a poset under inclusion, and in Section \ref{sec:embeddings} we give this poset structure for each simple factor. The subgroups are given up to conjugacy in the factor, except for type $D_7$ where for space reasons we list subgroups up to conjugacy in the full automorphism group of the factor. A Levi subgroup of type $D_7$ occurs for exceptional $G$ only when $G = E_8$, in which case $G$ induces the full automorphism group of this subgroup. In Section \ref{sec:conjugacy} we describe how the conjugacy classes in the other factors fuse under outer automorphisms.
\end{remark}

Note that when $L$ has a simple factor of exceptional type (including the case $L = G$), the image of $X$ in this exceptional group is irreducible, and therefore appears in existing work of the second author \cite{Tho1,Tho2,Tho3}. Since we give the poset of irreducible subgroups in each classical simple factor of each Levi subgroup of $G$, one can therefore construct the entire poset of $G$-cr connected subgroups from Theorem \ref{thm:main} and these references (see also Lemma \ref{lem:nonsimple_levi}). In particular, if the characteristic of the underlying field is zero or sufficiently large then every connected reductive subgroup of $G$ is $G$-cr (cf.\ \cite[Theorem 3.48]{MR2178661}), and thus all reductive subgroups of $G$ are known in this case.

We illustrate the above in Section \ref{sec:applying}, by constructing the poset of reducible $G$-completely reducible subgroups when $G$ has type $F_4$. We also give a series of lemmas regarding $G$-completely reducible subgroups, allowing us to derive much information on the subgroups which occur. The following summarises this information.

\begin{theorem} \label{thm:f4}
Let $G = F_{4}(K)$, where $K$ is algebraically closed of arbitrary characteristic, and let $X$ be a connected $G$-completely reducible subgroup of $G$. If $X$ is neither $G$-irreducible nor a torus, then the semisimple part of $X$ is conjugate to precisely one of the subgroups listed in Table \ref{tab:f4}, each of which is $G$-completely reducible.

For each semisimple subgroup $X$ appearing in Table \ref{tab:f4}, we give the action of $X$ on the adjoint module $L(G)$, the connected centraliser $C_{G}(X)^{\circ}$, as well as whether $X$ is separable in $G$ and whether $(G,X)$ is a reductive pair (see Section \ref{sec:action_sep_rp} for more details and definitions).
\end{theorem}

In addition to its illustrative purpose, Theorem \ref{thm:f4} also corrects some omissions in \cite[Corollary 5]{MR3075783}; more details are given in Remark \ref{rem:f4} on page \pageref{rem:f4}.

\section{Notation} \label{sec:nota}

We now present notation used throughout. Since this paper complements \cite{Tho3} we remain consistent with the notation there, which we recall for the convenience of the reader.

Throughout, algebraic groups are affine and defined over an algebraically closed field $K$. We define $p$ to be the characteristic of $K$ when this is positive, and set $p = \infty$ when $K$ has characteristic zero, so that expressions of the form `$p \ge n$' and `$p \neq n$' include the case that $K$ has characteristic zero. Subgroups are taken to be Zariski-closed, and representations are taken to be rational. Let $G$ be a simple algebraic group. We fix a maximal torus $T$ of $G$ and corresponding set of roots $\myPhi$. Let $\myPhi^{+}$ be a choice of positive roots, with corresponding simple roots $\myPi = \{ \alpha_1, \ldots, \alpha_l \}$ and fundamental dominant weights $\{ \lambda_1, \ldots, \lambda_l \}$. We use Bourbaki numbering \cite[p.\ 250]{MR0240238}, and sometimes use $a_1 a_2 \ldots a_l$ to denote a weight $a_1 \lambda_1 + a_2 \lambda_2 + \cdots + a_l \lambda_l$. When no ambiguity is possible, we write $V + W$ for the direct sum of $G$-modules $V$ and $W$, and $V^{*}$ denotes the dual module of $V$. We denote by $V_G(\lambda)$ (or just $\lambda$) the irreducible $G$-module of high weight $\lambda$. Similarly, the Weyl module and tilting module of high weight $\lambda$ will respectively be denoted $W(\lambda)$ and $T(\lambda)$. Given dominant weights $\mu_1, \mu_2, \ldots, \mu_n$, we define $T(\mu_1; \mu_2; \ldots ; \mu_n)$ to be the tilting module with the same composition factors as $W(\mu_1) + W(\mu_2) + \cdots + W(\mu_n)$, when this exists and is uniquely determined. For example, when $G$ is of type $A_1$ we use the notation $T(6;2)$ for $p \geq 5$ only. In this case $T(6;2) = 6 + 2$ when $p \geq 7$ and $T(6;2) = T(6)$ when $p = 5$. We use $L(G)$ to denote the adjoint module for $G$. Let $J = \{\alpha_{j_1}, \alpha_{j_2}, \ldots, \alpha_{j_r}\} \subseteq \myPi$ and define $\myPhi_J = \myPhi \cap \mathbb{Z}J$. Then the standard Levi subgroup corresponding to $J$ is $L_{j_1 j_2 \ldots j_r} = \left< T, U_\alpha \, : \, \alpha \in \myPhi_J \right>$. The notation $\bar{A}_{n}$ denotes a subgroup of $G$ of type $A_n$ which is generated by long root subgroups of $G$, and if $\myPhi$ admits multiple root lengths then $\tilde{A}_{n}$ denotes a subgroup of type $A_n$ which is generated by short root subgroups. This helps to distinguish between multiple subgroups of type $A_{n}$ which arise together. No ambiguity arises when considering simple subgroups of other types, and so we do not use such notation in these cases.

Now let $G = G_1 \ldots G_m$, where each $G_i$ is simple, and suppose $K$ has prime characteristic $p$. Let $F$ be the standard Frobenius endomorphism of $G$ acting on root groups $U_\alpha = \{ u_\alpha(c) \mid c \in K\}$ by $F(u_{\alpha}(c)) = u_{\alpha}(c^{p})$, and let $\rho \,:\, G \rightarrow \GL(V)$ be a representation. Then for a non-negative integer $r$, $V^{[r]}$ denotes the `twisted' module afforded by the representation $\rho^{[r]} \deq \rho \circ F^r$. 

We let $(V_1, \ldots, V_m)$ denote the $G$-module $V_1 \otimes \cdots \otimes V_m$, where $V_i$ is a $G_i$-module for each $i$, and we often use $0$ to denote the trivial module for $G$. Let $M_1, \ldots, M_k$ be $G$-modules and $n_1, \ldots, n_k$ be positive integers. Then $M_1^{n_1} / \ldots / M_k^{n_k}$ denotes a $G$-module having the same composition factors as $M_1^{n_1} + \dots + M_k^{n_k}$. Furthermore, $V = M_1 \mid \cdots \mid M_k$ denotes a $G$-module with socle series $V = V_{1} \supseteq V_{2} \supseteq \ldots \supseteq V_{k+1} = \{0\}$, with $V_{i}/V_{i+1} \cong M_{i}$ for $i = 1,\ldots,k$. The following notation is used in Section \ref{sec:tables}. Suppose $\mu_i$ is a dominant weight for each $G_i$. Then we define
\[ \myTheta(\mu_1; \ldots ;\mu_n) \deq 0 \mid ((V_{G_1}(\mu_1),0, \ldots, 0) + \cdots + (0, \ldots, 0, V_{G_n}(\mu_n)) ) \mid 0.\]
In other words, $\myTheta(\mu_1; \ldots ;\mu_n)$ has a unique minimal submodule, which is trivial, and a unique maximal submodule, with trivial quotient, and the quotient of the maximal submodule by the minimal submodule is a direct sum of irreducible modules, one for each simple factor of $G$.

\subsection{Diagonal subgroups} Many subgroups encountered here are \emph{diagonal}, in the following sense. If $H$ is a simply-connected simple group, then a diagonal subgroup of $H^{n} = H \times \cdots \times H$ is defined as the image of a morphism $H \to H^{n}$, $h \mapsto (\phi_1(h),\ldots,\phi_n(h))$, where each $\phi_i$ is a surjective endomorphism of $H$. A general semisimple group $Y$ is a central quotient of a direct product $H_{1}^{n_1} \times \cdots \times H_{s}^{n_s}$, for simply-connected simple groups $H_i$, and we define a diagonal subgroup to be the image of a morphism
\[ H_1 \times \cdots \times H_s \xrightarrow{\psi_1 \times \cdots \times \psi_s} H_{1}^{n_1} \times \cdots \times H_{s}^{n_s} \twoheadrightarrow Y, \]
where each $\psi_i \colon H_{i} \to H_{i}^{n_i}$ gives a diagonal subgroup in the previous sense.

Thus to specify such a diagonal subgroup, we must specify a surjective homomorphism from $H_i$ to each simple factor of $H_i^{n_i}$, for each $i$. Since $H_i$ is isomorphic to each factor of $H_{i}^{n_i}$, this is equivalent to specifying a surjective endomorphism of $H_i$ for each factor. By \cite[Section 1.15]{MR1490581}, such an endomorphism has the form $\alpha \theta F^{r}$ for an inner automorphism $\alpha$, a graph endomorphism $\theta$ and a power $F^{r}$ of a standard Frobenius endomorphism $F$. Since we only distinguish diagonal subgroups up to conjugacy, we can always take $\alpha$ to be trivial. Moreover a graph endomorphism induces an endomorphism on the weight lattice, and unless $H_i$ is of type $D_n$ ($n \geq 4$), the graph endomorphism is determined up to inner automorphism by the image of the weight $\lambda_1$; this includes the exceptional graph morphisms of $B_2$ and $F_4$ when $p = 2$ and of $G_2$ when $p=3$, which respectively send $\lambda_1$ to $2\lambda_2$, $2\lambda_4$ and $3\lambda_2$.

Our notation for diagonal subgroups is therefore as follows. When the simply-connected cover of a semisimple group $Y$ is $H^{n}$ for a single simply-connected group $H$, we denote a diagonal subgroup $X$ of $Y$ by
\begin{align*}
X \hookrightarrow Y \text{ via } (\mu_1^{[r_1]}, \ldots, \mu_k^{[r_k]}),
\end{align*}
where the $i$-th endomorphism is $\theta_i F^{r_i}$, with $\theta_i$ a (possibly trivial) graph endomorphism sending $\lambda_1$ to $\mu_i$, and $F$ a standard Frobenius morphism. In Section \ref{sec:tables} the group $Y$ is always clear and so we drop `$\hookrightarrow Y$' from the notation there. If $Y_1 Y_2 \ldots Y_s$ is a commuting product of groups of the form $Y$ above, then we concatenate the tuples $(\mu_1^{[r_1]}, \ldots, \mu_k^{[r_k]})$ to specify the diagonal subgroup, i.e.\ we write
\begin{align*}
X_1 X_2 \ldots X_s \hookrightarrow Y_1 Y_2 \ldots Y_s \text{ via } (\mu_{11}^{[r_1]}, \ldots, \mu_{1k_1}^{[r_{k_1}]}, \ldots, \mu_{sk_s}^{[r_{k_s}]}).
\end{align*}
Finally, if two or more of the $Y_{i}$ consist of factors of the same type, we must be careful to distinguish where each factor of $X_i$ is sent. We do this with subscripts $a$, $b$, etc. For example, if we write
\begin{align*}
X = A_{1}^{2} \hookrightarrow A_{1}^{4} = Y \text{ via } (1_a,1^{[r]}_a,1_b,1^{[s]}_b)
\end{align*}
then the subscript `$a$' in the first two coordinates means that the first factor $A_{1}$ of $X$ is a diagonal subgroup of the first two factors $A_{1}^{2}$ of $Y$, via $(1,1^{[r]})$, and the subscript `$b$' indicates that the second factor of $X$ is a diagonal subgroup of the last two factors $A_1^{2}$ of $Y$, via $(1,1^{[s]})$.

\section{Preliminaries} \label{sec:prelims}

We now give a series of results used in the proof of Theorem \ref{thm:main}. In Section \ref{sec:identifying} we give a criterion for identifying irreducible subgroups of simple algebraic groups of classical type. In Section \ref{sec:maxsubs} we enumerate the maximal connected subgroups of simple algebraic groups of types $A_6$ and $D_7$, which are used in deriving the poset structure of subgroups in Section \ref{sec:tables}. The maximal connected subgroups of other classical groups can be found in \cite[Lemma 3.3]{Tho3}.

\subsection{Identifying $L_{0}$-irreducible subgroups} \label{sec:identifying}

The following is central in our proof of Theorem \ref{thm:main}, as it provides a method of identifying $L_{0}$-irreducible subgroups when $L_{0}$ has classical type. The result is essentially well-known, but is usually stated in the literature as a one-way implication, see for example \cite[Lemma 2.2]{MR2043006}. We therefore provide a proof for completeness.

\begin{proposition} \label{prop:irredclass}
Let $G$ be a simple algebraic group of classical type over the algebraically closed field $K$, let $V = V_{G}(\lambda_1)$, and let $X$ be a subgroup of $G$. Then $X$ is $G$-irreducible if and only if one of the following holds:
\begin{enumerate}[label=\normalfont(\roman*)]
\item $G$ has type $A_n$ and $X$ acts irreducibly on $V$;
\item $G$ has type $B_n$, $C_n$ or $D_n$ and $V \downarrow X = V_1 \perp \ldots \perp V_k$ with the $V_i$ all non-degenerate, irreducible and pairwise inequivalent as $X$-modules;
\item $K$ has characteristic $2$, $G$ has type $D_{n}$ and $X$ fixes a non-singular vector $v \in V$, such that $X$ is $G_v$-irreducible in the point stabiliser $G_v$ and does not lie in a subgroup $D_{n-1}$ of $G_{v}$.
\end{enumerate}
\end{proposition}

\begin{remark}
In part (iii), the subgroup $G_v$ is simple of type $B_{n-1}$, and $X$ lies in a subgroup $D_{n-1}$ of $G_v$ if and only if the non-degenerate bilinear form on the $G_v$-module $\left<v\right>^{\perp}/\left<v\right>$ is the polarisation of an $X$-invariant quadratic form.
\end{remark}

\proof It is well-known that for $G = \SL(V)$, $\Sp(V)$ or $\SO(V)$, the parabolic subgroups of $G$ are respectively the stabilisers of flags of subspaces, flags of totally isotropic subspaces, and flags of totally singular subspaces \cite[Proposition 12.13]{MR2850737}. Thus if (i) holds then $X$ is $G$-irreducible. Similarly, if (ii) holds then every $X$-invariant subspace of $V$ is a direct sum of some of the factors $V_i$. Every such subspace being non-degenerate, it follows that $X$ is $G$-irreducible.

Next, for a contradiction suppose that (iii) holds and $X$ lies in a proper parabolic subgroup of $G$. Then $X$ stabilises a non-zero totally singular subspace $W$ of $V$. Note that $W \cap \left<v\right> = \{0\}$ since $v$ is non-singular, and so $W \cap \left<v\right>^{\perp}$ maps isomorphically onto its image in $\left<v\right>^{\perp}/\left<v\right>$. This image is a totally isotropic subspace of $\left<v\right>^{\perp}/\left<v\right>$, and since $X$ is $G_v$-irreducible we conclude that $W \cap \left<v\right>^{\perp} = \{0\}$. Since $\left<v\right>^{\perp}$ has codimension $1$ in $V$, this implies that $W$ has dimension at most $1$, hence exactly $1$ as $W$ is non-zero. So we may write $W = \left<w\right>$ where $(v,w) \neq 0$. Now $(W + \left<v\right>)^{\perp}$ is an $X$-submodule contained with codimension $1$ in $\left<v\right>^{\perp}$. Suppose that for some $c,d \in K$ we have $(av + bw,cv + dw) = 0$ for all $a,b \in K$. As $w$ is singular, this expands to $ac(v,v) + (ad + bc)(v,w) = 0$. Since $(v,w) \neq 0$, evaluating this at $a = 0$ and any $b \neq 0$ shows that $c = 0$. So the equation simplifies to $ad(v,w) = 0$ for all $a \in K$, showing that $d = 0$ also. Thus $(W + \left<v\right>)^{\perp} \cap (W + \left<v\right>) = \{0\}$, so $(W + \left<v\right>)^{\perp}$ is non-degenerate. This contradicts the assumption that $X$ does not lie in a subgroup of type $D_{n-1}$ in $G_v$, and we conclude that $X$ is $G$-irreducible.

Conversely, suppose $X$ is $G$-irreducible and that neither (i) nor (iii) hold. If $K$ has characteristic $2$, $G = D_n$ and $X$ fixes a non-singular vector $v$ on $V$, then since (iii) does not hold, $X$ lies either in a subgroup $D_{n-1}$ of $G$, or in a parabolic subgroup of $G_v$; each of these contradicts the hypothesis that $X$ is $G$-irreducible. So $X$ does not fix any non-singular vector of $V$. Thus the bilinear form is non-degenerate on each irreducible $X$-submodule, so each irreducible $X$-submodule has an orthogonal complement, and $V \downarrow X$ is completely reducible. Let $U$ be an irreducible $X$-submodule of $V$, and suppose that $W \subseteq U^{\perp}$ is an irreducible $X$-submodule with $U \cong W$. Since $U$ and $W$ are irreducible, their $X$-invariant non-degenerate bilinear or quadratic forms are uniquely determined up to a constant, hence there exists an isometric $X$-module isomorphism $\phi \, : \, U \to W$. For a scalar $\lambda \in K$, consider the $X$-invariant subspace $\{ u + \lambda\phi(u) \, : \, u \in U \}$. Since $U$ and $W$ are orthogonal, we have $(u + \lambda\phi(u),u' + \lambda\phi(u')) = (1 + \lambda^{2})(u,u')$. In particular, taking $\lambda$ to be a square root of $-1$ in $K$, we obtain a totally isotropic submodule. This contradiction shows that $U^{\perp}$ does not contain any submodule isomorphic to $U$, and thus (ii) holds. \qed

\begin{lemma} \label{lem:nonsimple_levi}
Let $L$ be a Levi subgroup of $G$. Then a subgroup $X$ of $L$ is $L$-irreducible if and only if the image of $X$ under projection to each simple factor of $L/Z(L)$ is irreducible.
\end{lemma}

\proof This follows because the parabolic subgroups of $L$ are precisely the products of $Z(L)$ with parabolic subgroups of the simple factors of $L$, cf.\ \cite[Lemmas 2.11, 2.12]{MR2178661}. \qed

\subsection{Maximal subgroups of classical groups of low rank} \label{sec:maxsubs}

Proposition \ref{prop:irredclass} lets us identify irreducible subgroups of classical groups. Since we are also interested in the poset of $G$-cr subgroups, it will be of use for us to know which subgroups arising in Proposition \ref{prop:irredclass} are maximal among connected subgroups. This allows us to work recursively through chains of maximal connected subgroups. For many groups of low rank, this information has been calculated in \cite[Lemma 3.3]{Tho3}; we additionally require the maximal subgroups of simple algebraic groups of types $A_6$ and $D_7$. We begin with a lemma on fusion of conjugacy classes under outer automorphisms.
\begin{lemma} \label{lem:go_so_classes}
Let $G$ be a group with a normal subgroup $S$ of finite index, and let $X$ be a subgroup of $S$. Then the conjugacy class $X^{G}$ splits into $|G : S|$ classes of subgroups of $S$ if and only if $N_{G}(X) \le S$. In particular, if $|G:S|$ is prime then $S$ is transitive on $X^{G}$ if and only if $X$ is normalised by an element of $G \setminus S$, otherwise $X^{G}$ splits into $|G:S|$ classes of subgroups of $S$.
\end{lemma}

\proof Consider the transitive action of $G$ on the set $X^{G}/S$ of $S$-conjugacy classes of $G$-conjugates of $X$. A point stabiliser is conjugate to $S N_{G}(X)$, and therefore $|X^{G}/S| = |G:SN_{G}(X)|$, which equals $|G:S|$ if and only if $N_{G}(X) \le S$. \qed

Our classification is based on \cite{MR1650328}, which separates maximal subgroups of classical groups into classes, of which we require the following. Write $G = Cl(V)$ to denote that $G$ is one of the groups $\SL(V)$, $\SO(V)$ or $\Sp(V)$.

{\bf Class} $\mathcal{C}_{1}$: Subspace stabilisers. Here $H \in \mathcal{C}_{1}$ if $H = \Stab_{G}(W)$ where $W$ is either a non-degenerate subspace of $V$ ($V$ is given the zero form if $G = \SL(V)$), or $(G,p) = (\SO(V),2)$ and $W$ is a non-singular subspace of dimension $1$.

{\bf Class} $\mathcal{C}_{4}$: Tensor product subgroups. Suppose that $V = V_{1} \otimes V_{2}$ with $\dim V_{i} > 1$ for each $i$. Then $H \in \mathcal{C}_{4}$ if $H = Cl(V_1) \circ Cl(V_2)$, acting on $V$ by $(g_1,g_2)(v_1 \otimes v_2) \deq (g_1 v_1) \otimes (g_2 v_2)$. The subgroups occurring here have one of the following forms:
\begin{align*}
\SL \otimes \SL &< \SL, & \Sp \otimes \SO &< \Sp\ (p \neq 2),\\
\Sp \otimes \Sp &< \SO, & \SO \otimes \SO &< \SO\ (p \neq 2).
\end{align*}

The following is immediate from \cite[Theorem 1]{MR1650328} (cf.\ also \cite[Theorem 3.2]{Tho3}). Here a restricted irreducible module is a module whose highest weight, when expressed as a sum of fundamental dominant weights, has coefficients strictly less than $p$.

\begin{lemma} \label{lem:maxclassical}
Let $G = Cl(V)$ be a classical simple algebraic group. Suppose that $M$ is a reductive, maximal connected subgroup of $G$. Then one of the following holds:
\begin{enumerate}[label=\normalfont(\roman*)]
	\item $M$ belongs to $\mathcal{C}_{1}$;
	\item $M$ belongs to $\mathcal{C}_{4}$;
	\item $M$ is a simple algebraic group and $V \downarrow M$ is irreducible and restricted.
\end{enumerate}
\end{lemma}

\begin{proposition} \label{prop:max_a6_d7}
Let $G$ be a classical simple algebraic group of type $A_6$ or $D_7$. Then Table \ref{tab:maxclassical} gives all $G$-conjugacy classes of reductive, maximal connected subgroups of $G$.
\end{proposition}

\proof For $G$ of type $A_6$, the class $\mathcal{C}_{1}$ contains no $G$-irreducible subgroups, and the class $\mathcal{C}_{4}$ is empty since the natural module is $7$-dimensional. Thus a maximal subgroup $M$ falls into case (iii) of Lemma \ref{lem:maxclassical}. The only simple algebraic groups with a $7$-dimensional irreducible module are those of type $A_6$, and of type $B_3$ and $G_2$ when $p \neq 2$. The module is then unique up to taking a Frobenius twist, and so $G$ contains a unique subgroup of each of these types up to conjugacy. Moreover a subgroup of type $B_3$ contains a subgroup of type $G_2$, and it follows that every subgroup of $G$ of type $G_2$ is contained in a subgroup of type $B_3$.

For $G$ of type $D_7$, the $G$-irreducible subgroups in class $\mathcal{C}_{1}$ are of type $B_6$, $A_1 B_5$, $\bar{A}_1^2 D_5$, $B_2 B_4$, $\bar{A}_3 D_4$ and $B_3^2$, respectively corresponding to orthogonal decompositions of the $14$-dimensional module into spaces of dimension $i$ and $14 - i$ for $i \in \{1,3,4,5,6,7\}$ when $p \neq 2$. When $p = 2$, the subgroups $A_1 B_5$, $B_2 B_4$ and $B_3^2$ each fix a non-singular vector, and are therefore contained in a subgroup $B_6$.

The prime decomposition $14 = 2 \times 7$, and the fact that there does not exist a non-degenerate orthogonal module of dimension $2$ or symplectic module of dimension $7$, shows that the class $\mathcal{C}_{4}$ is empty. It remains to consider the subgroups in case (iii) of Lemma \ref{lem:maxclassical}. By \cite{MR1901354} the restricted, irreducible, $14$-dimensional modules for a simple algebraic group are as follows:
$V_{B_2}(20)$ $(p \neq 2,5)$, 
$V_{G_2}(01)$ $(p \neq 3)$,
$V_{A_3}(101)$ $(p = 2)$,
$V_{B_3}(010)$ $(p = 2)$,
$V_{C_3}(010)$ $(p \neq 3)$, and
$V_{C_3}(001)$ $(p \neq 2)$. 
Using \cite[Lemma 79]{steinberglectureson}, it is straightforward to calculate that if $p \neq 2$, then all but the last of these modules supports a non-degenerate quadratic form, and therefore gives rise to a subgroup of $G$ (the final module is symplectic). If $p = 2$ then inspecting \cite[Table 1, pp.\ 282--283]{MR888704}, in particular numbers ${\rm IV}_{7}$, ${\rm IV}_{8}$ and ${\rm S}_{7}$ there, shows that the relevant modules here still support a non-degenerate quadratic form. Note that, since $V_{B_3}(010)$ is a direct summand of $\bigwedge^{2}(V_{B_3}(100))$, the representation factors through the morphism $B_{3} \to C_{3}$, so the image in $D_{7}$ is a subgroup of type $C_{3}$, and is thus conjugate to the subgroup given by $V_{C_3}(010)$. When $p = 2$, a subgroup $C_3$ contains a subgroup $A_3$ acting via $V_{A_3}(101)$ and a subgroup $G_2$ acting via $V_{G_2}(01)$, and so these latter subgroups are not maximal in $G$. Thus when $p = 2$ the only possible maximal subgroups as in part (iii) of Lemma \ref{lem:maxclassical} are of type $C_3$. When $p \neq 2$, the possible subgroups (of type $B_2$, $G_2$ and $C_3$) do not contain irreducible copies of one another, and so each gives rise to a maximal subgroup of $G$.

There is a unique subgroup class of each of these three types in $\GO_{14}(K)$, and each has trivial group of algebraic outer automorphisms. Thus if $X$ is such a subgroup, then either $N_{\GO_{14}(K)}(X) \le \SO_{14}(K)$ or $X$ centralises a non-trivial element of $\GO_{14}(K) \setminus \SO_{14}(K)$. The latter does not occur, since the centraliser of such an element is the stabiliser of an orthogonal decomposition of the natural module, and $X$ does not preserve such a decomposition. Thus by Lemma \ref{lem:go_so_classes}, $X$ gives rise to two classes of subgroups of $\SO_{14}(K)$, hence to two classes of subgroups of $G$. \qed

{ \small

\begin{table}
\centering
\caption{Maximal connected subgroups of certain classical groups.}
\label{tab:maxclassical}
\begin{tabular}{llll} \hline
$G$ & Maximal subgroup $M$ & $V_{G}(\lambda_1) \downarrow M$ & Comments \\ \hline
$A_6$ & $B_3$ $(p \neq 2)$ & $100$ & \\ \hline
$D_7$ & $B_6$ & $T(\lambda_1;0)$ & \\
& $A_1 B_5$ $(p \neq 2)$ & $(2,0) + (0,\lambda_1)$ & \\
& $\bar{A}_1^2 D_5$ & $(1,1,0) + (0,0,\lambda_1)$ & \\
& $B_2 B_4$ $(p \neq 2)$ & $(10,0) + (0,\lambda_1)$ & \\
& $\bar{A}_3 D_4$ & $(010,0) + (0,\lambda_1)$ & \\
& $B_3^2$ $(p \neq 2)$ & $(100,0) + (0,100)$ & \\
& $C_3$ $(p \neq 3)$ (2 classes) & $010$ & classes permuted by graph aut. \\
& $B_2$ $(p \neq 2, 5)$ (2 classes) & $20$ & classes permuted by graph aut. \\
& $G_2$ $(p \geq 5)$ (2 classes) & $01$ & classes permuted by graph aut. \\ \hline
\end{tabular}
\end{table}
}

\section{Proof of the main theorem} \label{sec:pfthm}

In this section we describe the method used for proving Theorem \ref{thm:main}, i.e.\ for determining the poset of $L_{0}$-irreducible connected subgroups of all classical simple factors $L_{0}$ of a Levi subgroup of an exceptional algebraic group $G$. It is very similar to the strategy described in \cite[Section 4]{Tho3}, and employs a system of ID numbers for the different conjugacy classes of $L_0$-irreducible connected subgroups. The method is identical for each factor $L_{0}$ and so we describe the general method, and give explicit details only for $L_{0}$ of type $A_3$. 

Fix the factor $L_{0}$. Firstly, we find all of the reductive, maximal connected subgroups of $L_{0}$, given by Lemma \ref{lem:maxclassical} and \cite[Lemma~3.3]{Tho3}; these subgroups are automatically $L_{0}$-irreducible. Let $M$ be the first reductive, maximal connected subgroup listed. We now find the $L_{0}$-irreducible maximal connected subgroups $M_1$ of $M$. To do this, we again find the reductive, maximal connected subgroups of $M$, using \cite[Lemma~3.3]{Tho3} for classical factors of $M$ and \cite[Corollary~12]{MR898346} for factors of type $G_2$. Using Proposition \ref{prop:irredclass} and Lemma \ref{lem:nonsimple_levi}, we check whether these maximal connected subgroups are $L_{0}$-irreducible. We now repeat this process for each $L_{0}$-irreducible maximal connected subgroup of $M$. Continuing in this manner yields all $L_{0}$-irreducible subgroups of $L_{0}$. However, we must be careful to avoid repeats. Thus at each step, we check to see if each subgroup $X$ arising has been previously found. If not, we assign the subgroup a new unique ID number, and if it has been found before, we use its existing ID number. To check for repeats we consider the action of $X$ on $V_{L_{0}}(\lambda_1)$. This determines the $L_{0}$-class of $X$, unless perhaps if $L_{0}$ has type $D_n$, when the image of $X$ is determined up to $\Aut(D_n)$-conjugacy and we use Lemma \ref{lem:go_so_classes} to decide how these split in $L_{0}$. At this stage, we now know all of the connected overgroups of each $L_{0}$-irreducible subgroup, and thus understand the poset structure of the irreducible connected subgroups of $L_{0}$.

The information obtained in this way is given in Tables \ref{A2tab}--\ref{D7tab}. The notation used there is explained in Section~\ref{sec:nota}, and we explain how to read the tables at beginning of Section~\ref{sec:tables}. In particular, when an ID number is given in italics, it means the corresponding subgroup has appeared elsewhere in the table. 

To save space, we deviate from the above slightly when considering diagonal subgroups. Suppose that $X$ is an $L_{0}$-irreducible connected subgroup of the form $A^n B$ for some $n \geq 2$ with $A$ and $B$ of different types. Then $X$ has maximal diagonal subgroups of the form $A^{n-1} B$, and these have subgroups of the form $A^{n-2}B$, and so on. We list all such diagonal subgroups of $X$ at once. This significantly reduces the size of the tables, without missing any $L_{0}$-irreducible subgroups. The cost is that some additional combinatorial work is required to recover the poset of overgroups of certain diagonal subgroups. For readability, in some cases we move diagonal subgroups to a supplementary table in Section \ref{tabs:diags}.

There is a subtlety concerning the position of a subgroup in the poset of irreducible subgroups. For instance, consider $L_{0} = \SO_{8}$, of type $D_4$. When $p \neq 2$ this contains a maximal subgroup $\SO_5 \SO_3$. When $p=2$ such a subgroup still exists and is $\SO_8$-irreducible, but is now contained in $\SO_7$. Thus we use the same ID number $n$ for both groups but write $n$a when $p \neq 2$ and $n$b when $p = 2$. This allows us to discuss both cases together. In the tables, if we first arrive at the subgroup labelled $n$b, we postpone listing its subgroups until reaching the subgroup labelled $n$a. This generalises to any situation where `the same' subgroup $X$ occurs in different places in the poset depending on $p$. In particular, whenever we use this notation the subgroups $n$a and $n$b have identical weights (with multiplicities) on $V_{L_0}(\lambda_1)$.

\subsection{The irreducible subgroups of $L_{0} = A_3$}

The reductive, maximal connected subgroups of $A_3$ are $B_2$ and $A_1^2$ $(p \neq 2)$, respectively acting on $V = V_{A_3}(\lambda_1)$ via $01$ and $(1,1)$. These are $L_{0}$-irreducible by Proposition \ref{prop:irredclass}. The module $(1,1)$ yields an irreducible $A_1^2$ subgroup of $A_3$. When $p = 2$ we will shortly find this subgroup inside $B_2$, and this explains why we write $2$a for the ID number of $A_1^2$ $(p \neq 2)$ in Table \ref{A3tab}. Next, we consider the reductive, maximal connected subgroups of $B_2 = A_3(\#1)$. These are $A_1^2$, $A_1^2$ $(p=2)$ and $A_1$ $(p \geq 5)$ acting as $(1,0) + (0,1)$, $(1,1)$ and $3$ on $V$, respectively. The former $A_1^2$ subgroup is not $L_{0}$-irreducible by Proposition \ref{prop:irredclass}, since it acts reducibly on $V$. The second $A_1^2$ subgroup is $A_3$-irreducible and given ID number $2$b because of the explanation above. The subgroup $A_1$ $(p \geq 5)$ is also $A_3$-irreducible and given ID number $3$. We now need to consider the $A_3$-irreducible connected subgroups of $A_3(\#2\textup{b}) = A_1^2$ $(p=2)$ and $A_3(\#3) = A_1$. As discussed in Section \ref{sec:tables}, for presentation reasons we postpone studying the irreducible subgroups of $A_3(\#2\text{b})$ and instead consider the subgroups of $A_3(\#2\text{a})$ and $A_3(\#2\text{b})$ together later. Since $A_1$ contains no proper reductive, maximal connected subgroups, there is nothing more to list so we turn immediately to $A_3(\#2) = A_1^2$ (this means the class of subgroups acting via $(1,1)$ with no restriction on the characteristic). The reductive, maximal connected subgroups are diagonal $A_1$ subgroups. These are embedded via $(1^{[r]},1^{[s]})$ $(rs=0)$. These then act on $V$ via $1^{[r]} \otimes 1^{[s]}$. Therefore, they are $A_3$-irreducible if and only if $r \neq s$. Furthermore, there is an element of $N_{A_3}(A_1^2) = (A_1^2).2$ swapping the two $A_1$ factors. Therefore, up to $A_3$-conjugacy we have the $L_{0}$-irreducible subgroups $A_1 \hookrightarrow A_1^2$ via $(1,1^{[r]})$ $(r \neq 0)$. These we assign the ID number $4$.

\section{Tables of subgroups} \label{sec:tables}

We now present the embeddings referred to in Theorem \ref{thm:main}, and begin by describing how to read the tables. There is one table for each type of simple classical factor $L_{0}$ of a Levi subgroup of an exceptional algebraic group, except for $L_{0}$ of type $A_1$ where no table is required. Since we are considering only connected subgroups, there in no harm in replacing $L_{0}$ by an isogenous group (finite central extension or quotient), which we do without further comment.

\subsection{Identification numbers} In order to discuss the poset of subgroups occurring, we assign each class of subgroups an ID number. Each subgroup listed denotes either a unique $L_{0}$-conjugacy class of subgroups, or an infinite collection of classes parametrised by some Frobenius twists. In Table \ref{D7tab} subgroups are instead given up to $\Aut(D_7)$-conjugacy for space reasons. We write $L_{0}(\#n)$ (or just $n$ when $L_{0}$ is clear) to refer to the class of $L_{0}$-irreducible subgroups with ID number $a$. We set $L_{0}(\#0)$ to be $L_{0}$ itself.

If $L_{0}(\#n)$ refers to a collection of classes depending on Frobenius twists with exponents $r_1,\ldots,r_k$ and we wish to refer to a proper subset of these classes, we write $L_{0}(\#n^{\{r_1, \ldots, r_k\}})$ to specify the powers in the field twists.

Let us give an illustration. Consider the line $A_5(\#8)$ in Table \ref{A5tab}, giving diagonal subgroups $A_1 \hookrightarrow A_1 A_1$ via $(1^{[r]},1^{[s]})$ $(p \neq 2; rs = 0; r \neq s$). Then $A_5(\#8^{\{r,0\}})$ refers to the conjugacy classes with $s=0$ and the notation $A_5(\#{8}^{\{1,0\}})$ refers to the single conjugacy class $A_1 \hookrightarrow A_1 A_1$ via $(1^{[1]},1)$.

We sometimes use a shorthand for $n^{\{r_1, \ldots, r_k\}}$. The notation $\mathit{n^{\{\underline{0}\}}}$ simply means that each $r_i$ is equal to $0$. The notation $n^{\{\delta_j\}}$ means $r_j = 1$ and $r_i = 0$ for all $i \neq j$.

\subsection{Reading the tables} Each table is divided into sections, separated by horizontal lines. We first list the maximal $L_{0}$-irreducible connected subgroups, and in each subsequent section we list the $L_{0}$-irreducible maximal connected subgroups contained in a (previously listed) subgroup, say $X$, as well as all $L_{0}$-irreducible diagonal subgroups of $X$. The section heading gives the type and ID number of $X$, as well as any restrictions on the characteristic $p$ required for $X$ to exist and be $L_{0}$-irreducible. There is one more piece of information in the heading, which is the ``$M_i =$'' for $i = 0,1, \dots$ (where $M_0$ is simply written as $M$). This is intended to make it easier for the reader to follow the tables, and is explained below. 

Other than section headings, the information given in a row depends on whether the subgroup in question is a diagonal subgroup of its immediate overgroup $X$. If $Y$ is a diagonal subgroup of $X$, then the second column contains a statement ``$Y$ via \dots'' giving the embedding of $Y$ into $X$, as well as any restrictions on the characteristic or Frobenius twists in the embedding. If $Y$ is not diagonal then the second column gives the isomorphism type of $Y$ and restrictions on the characteristic, while the third gives the restriction $V_{L_{0}} \downarrow Y$, where $V_{L_0} = V_{L_{0}}(\lambda_1)$. For a diagonal subgroup $Y$ of $X$, it is straightforward to work out $V_{L_{0}} \downarrow Y$ from $V_{L_{0}} \downarrow X$. Restrictions on $p$ are inherited by subgroups, so for example, if $X$ exists and is $L_{0}$-irreducible only for $p \neq 2$, then we write $X$ $(p \neq 2)$ in the heading and any subgroup $Y$ of $X$ implicitly inherits this restriction on $p$. However, when we consider the subgroups of $Y$ we explicitly repeat any inherited restriction on $p$.

In Tables \ref{D6tab} and \ref{D7tab}, for readability we have moved large collections of diagonal subgroups to supplementary tables in Section \ref{tabs:diags}. This is marked by an entry ``See Table $x$'' where $x$ is the relevant table number.

We now explain how to determine the poset of $L_{0}$-irreducible connected subgroups using Tables \ref{A2tab}--\ref{D7tab}. Each table starts with the maximal $L_{0}$-irreducible connected subgroups, with their ID number listed in the first column. A pair of horizontal lines indicates the end of this list. Next we write ``In $M = H_1 H_2 \dots $ ($L_{0}(\#n_1)$)'' where $n_1$ is the identification number for the first reductive, maximal connected subgroup of type $H_1 H_2 \dots $, as well as any restrictions on the characteristic $p$. We then list the $L_{0}$-irreducible maximal connected subgroups of $M$, as well as any diagonal connected subgroups of $M$ (not just the maximal ones). Recall that we will not explicitly consider the proper subgroups of any diagonal connected subgroups, as discussed in Section \ref{sec:pfthm}. A horizontal line then indicates the end of this list. The next row then gives a heading ``In $M_1 = X_1 X_2 \dots$ ($L_{0}(\#n_2)$)'', where $X_1 X_2 \dots$ is the first $L_{0}$-irreducible maximal connected subgroup of $M$. The subscript in ``$M_1$'' indicates that we are now listing subgroups of a maximal connected subgroup of a maximal connected subgroup. We then repeat the process, listing the $L_{0}$-irreducible maximal connected subgroups of $M_1$ and all diagonal subgroups of $M_1$. The next heading could be ``In $M_2 = Y_1 Y_2 \dots $ ($L_{0}(\#n_3)$)'', where $Y_1 Y_2 \dots$ is a maximal connected subgroup of $M_1$ or it could be ``In $M_1 = Z_1 Z_2 \dots $ ($L_{0}(\#n_3)$)'', where $Z_1 Z_2 \dots$ is the second maximal connected subgroup of $M$. This will depend on whether $M_1$ has any proper $L_{0}$-irreducible connected subgroups. Once all $L_{0}$-irreducible connected subgroups of $M$ have been listed in this way, a pair of horizontal lines indicate the end of the subgroups contained in the first reductive, maximal connected subgroup of $G$. The next heading will be ``In $M = K_1 K_2 \dots$ ($L_{0}(\#n_4)$)'' and we repeat the process again for the second reductive, maximal connected subgroup $K_1 K_2 \dots$ of $L_{0}$. We iterate this process until we have considered all of the $L_{0}$-irreducible subgroups contained in all reductive, maximal connected subgroups.

To avoid redundancy, whenever some subgroup $X$ occurs as a maximal connected subgroup of two or more non-conjugate connected subgroups of $L_{0}$, in each occurrence of $X$ after the first we give the ID number in italics, and do not re-list the proper subgroups of $X$, nor do we repeat the restriction $V_{L_{0}} \downarrow X$.

There is another important scenario where we do not immediately list the subgroups contained in $X$. In this case the ID number will be $n$b and the subgroup $X$ is defined for all $p \not \in \{k_1, \dots, k_t\}$, but is only maximal in the subgroup currently being considered when the characteristic is some particular prime $l \not \in \{k_1, \ldots, k_t\}$. At some point later in the table, $X$ will be defined for all $p \not \in \{k_1, \dots, k_t\} \cup \{l\}$ and given ID number $n$a, with a single exception when $X = D_7(\#218)$. In this latter case, $X$ is a maximal subgroup of $D_7(\#183)$ when $p=2$ (where it gets the label $218$b) but $X$ is defined later in the table for all $p$ when it occurs as a maximal subgroup of $D_7(\#176)$ and simply given ID number $218$. When the subgroup with ID number $n$a (or $218$ in the exception) occurs, we then consider the subgroups of $X$ for all $p \not \in \{k_1, \dots, k_t\}$ together. There are instances where the subgroup $X$ occurs again in the table, in which case we simply write $\mathit{n}$ regardless of any characteristic restrictions.

\subsection{Tables of embeddings for Theorem \ref{thm:main}} \label{sec:embeddings} \leavevmode

{ \small 

\begin{longtable}{lll}
\caption{Irreducible subgroups of $L_{0} = A_2$. \label{A2tab}} \\
\hline ID & Subgroup $X$ & $V_{L_{0}} \downarrow X$ \\
\hline 1 & $A_1$ $(p \neq 2)$ & $2$ \\ \hline \hline \endfirsthead
\end{longtable}

\begin{longtable}{lll}
\caption{Irreducible subgroups of $L_{0} = A_3$. \label{A3tab}} \\
\hline ID & Subgroup $X$ & $V_{L_{0}} \downarrow X$ \\ \hline
1 & $B_2$ & $01$ \\
2a & $A_1^2$ $(p \neq 2)$ & $(1,1)$  \\ \hline \hline \endfirsthead
\multicolumn{3}{l}{In $M = B_2$ ($A_3(\#1)$)} \\
2b & $A_1^2$ $(p = 2)$ & $(1,1)$ \\
3 & $A_1$ $(p \geq 5)$ & $3$ \\ \hline \hline
\multicolumn{3}{l}{In $M = A_1^2$ ($A_3(\#2)$)} \\
4 & \multicolumn{2}{l}{$A_1$ via $(1,1^{[r]})$ $(r \neq 0)$} \\ \hline \hline
\end{longtable}

\begin{longtable}{lll}
\caption{Irreducible subgroups of $L_{0} = A_4$. \label{A4tab}} \\
\hline ID & Subgroup $X$ & $V_{L_{0}} \downarrow X$ \\
\hline 1 & $B_2$ $(p \neq 2)$ & $10$ \\ \hline \hline
\multicolumn{3}{l}{In $M = B_2$ $(p \neq 2)$ ($A_4(\#1)$)} \\
2 & $A_1$ $(p \geq 5)$ & $4$ \\ \hline \hline \endfirsthead
\end{longtable}

\begin{longtable}{lll}
\caption{Irreducible subgroups of $L_{0} = A_5$. \label{A5tab}} \\ \hline
ID & Subgroup $X$ & $V_{L_{0}} \downarrow X$ \\ \hline
1 & $C_3$ & $100$ \\
2 & $A_1 A_2$ & $(1,10)$ \\
3a & $A_3$ $(p \neq 2)$ & $010$ \\
4 & $A_2$ $(p \neq 2)$ & $20$ \\ \endfirsthead \hline \hline
\multicolumn{3}{l}{In $M = C_3$ ($A_5(\#1)$)} \\
5 & $A_1 A_1$ $(p \neq 2)$ & $(1,2)$  \\
3b & $A_3$ $(p=2)$ & $010$  \\
6 & $G_2$ $(p=2)$ & $10$  \\
7 & $A_1$ $(p \geq 7)$ & $5$  \\ \hline
\multicolumn{3}{l}{In $M_1 = A_1 A_1$ $(p \neq 2)$ ($A_5(\#5)$)} \\
8 & \multicolumn{2}{l}{$A_1$ via $(1^{[r]},1^{[s]})$ $(rs=0)$} \\ \hline \hline
\multicolumn{3}{l}{In $M = A_1 A_2$ ($A_5(\#2)$)} \\
\it{5} & $A_1 A_1$ $(p \neq 2)$ &  \\ \hline \hline
\end{longtable}

\begin{longtable}{lll}
\caption{Irreducible subgroups of $L_{0} = A_6$. \label{A6tab}} \\
\hline ID &Subgroup $X$ & $V_{L_{0}} \downarrow X$ \\ \hline
1 & $B_3$ $(p \neq 2)$ & $100$ \\ \hline \hline
\multicolumn{3}{l}{In $M = B_3$ $(p \neq 2)$ ($A_6(\#1)$)} \\
2 & $G_2$ & $10$ \\ \hline 
\multicolumn{3}{l}{In $M_1 = G_2$ $(p \neq 2)$ ($A_6(\#2)$)} \\
3 & $A_1$ $(p \geq 7)$ & $6$  \\
4 & $A_2$ $(p =3)$ & $11$  \\ \hline \hline \endfirsthead
\end{longtable}

\begin{longtable}{lll}
\caption{Irreducible subgroups of $L_{0} = A_7$. \label{A7tab}} \\ \hline
ID & Subgroup $X$ & $V_{L_{0}} \downarrow X$ \\ \hline
1 & $C_4$ & $\lambda_1$ \\
2a & $D_4$ $(p \neq 2)$ & $\lambda_1$ \\
3 & $A_1 A_3$ & $(1,100)$ \\ \hline \hline
\endfirsthead
\multicolumn{3}{l}{In $M = C_4$ ($A_7(\#1)$)} \\
4a & $A_1^3$ $(p \neq 2)$ & $(1,1,1)$  \\
5 & $A_1$ $(p \geq 11)$ & $7$  \\
2b & $D_4$ $(p=2)$ & $\lambda_1$  \\ \hline
\multicolumn{3}{l}{In $M_1 = A_1^3$ $(p \neq 2)$ ($A_7(\#4)$)} \\
6 & \multicolumn{2}{l}{$A_1 A_1$ via $(1_a,1_a^{[r]},1_b)$ $(r \neq 0)$} \\
7 & \multicolumn{2}{l}{$A_1$ via $(1,1^{[r]},1^{[s]})$ $(0 < r < s)$} \\ \hline \hline
\multicolumn{3}{l}{In $M = D_4$ ($A_7(\#2)$)} \\
8 & $B_3$ & $001$ \\
9a & $A_1 B_2$ $(p \neq 2)$ & $(1,01)$ \\
10 & $A_2$ $(p \neq 3)$ & $11$ \\ \hline
\multicolumn{3}{l}{In $M_1 = B_3$ ($A_7(\#8)$)} \\
9b & $A_1 B_2$ $(p = 2)$ & $(1,01)$ \\ \hline
\multicolumn{3}{l}{In $M_1 = A_1 B_2$ ($A_7(\#9)$)} \\
4b & $A_1^3$ $(p = 2)$ & $(1,1,1)$  \\ \hline \hline
\multicolumn{3}{l}{In $M = A_1 A_3$ ($A_7(\#3)$)} \\
\it{9} & $A_1 B_2$ &  \\ \hline \hline
\end{longtable}

\begin{longtable}{lll}
\caption{Irreducible subgroups of $L_{0} = C_2$. \label{C2tab}} \\ \hline
ID & Subgroup $X$ & $V_{L_{0}} \downarrow X$ \\ \hline
\ctwo{1} & $\bar{A}_1^2$ & $(1,0) + (0,1)$ \\
\ctwo{2} & $\tilde{A}_1^2$ $(p=2)$ & $(1,1)$ \\
\ctwo{3}a & $A_1$ $(p \geq 5)$ & $3$ \\ \hline \hline
\endfirsthead
\multicolumn{3}{l}{In $M = \bar{A}_1^2$ $(C_2(\#\ctwo{1}))$} \\
\ctwo{4} & \multicolumn{2}{l}{ $A_1$ via $(1,1^{[r]})$ $(r \neq 0)$   } \\ \hline \hline
\multicolumn{3}{l}{In $M = \tilde{A}_1^2$ $(p=2)$ $(C_2(\#\ctwo{2}))$} \\
\ctwo{3}b & \multicolumn{2}{l}{ $A_1$ via $(1,1^{[1]})$ } \\ 
\ctwo{5} & \multicolumn{2}{l}{ $A_1$ via $(1,1^{[r]})$ $(r > 1)$   } \\ \hline \hline
\end{longtable}

\begin{longtable}[c]{lll}
\caption{Irreducible subgroups of $L_{0} = C_3$. \label{C3tab}} \\ \hline
ID & Subgroup $X$ & $V_{L_{0}} \downarrow X$ \\ \hline
\cthree{1} & $\bar{A}_1 C_2$ & $(1,0) + (0,10)$ \\
\cthree{2} & $A_1 A_1$ $(p \neq 2)$ & $(2,1)$ \\

\cthree{3} & $\tilde{A}_3$ $(p=2)$ & $010$ \\
\endfirsthead
\cthree{4} & $G_2$ $(p=2)$ & $10$ \\ 

\cthree{5} & $A_1$ $(p \geq 7)$ & $5$ \\
\hline \hline
\multicolumn{3}{l}{In $M = \bar{A}_1 C_2$ $(C_3(\#\cthree{1}))$} \\

\cthree{6} & $\bar{A}_1 \tilde{A}_1^2$ $(p = 2)$ & $(1,0,0) + (0,1,1)$ \\

\cthree{7} & $\bar{A}_1 A_1$ $(p \geq 5)$ & $(1,0) + (0,3)$ \\

\cthree{8} & $\bar{A}_1^3$ & $(1,0,0) + (0,1,0) + (0,0,1)$ \\
\hline
\multicolumn{3}{l}{In $M_1 = \bar{A}_1 \tilde{A}_1^2$ $(p=2)$ $(C_3(\#\cthree{6}))$} \\
\cthree{9} & \multicolumn{2}{l}{ $A_1 \tilde{A}_1$ via $(1_a^{[r]},1_a^{[s]},1_b)$ ($rs=0$) } \\
\cthree{10} & \multicolumn{2}{l}{ $\bar{A}_1 A_1$ via $(1_a,1_b,1_b^{[r]})$ ($r \neq 0$) } \\
\cthree{11} & \multicolumn{2}{l}{ $A_1$ via $(1^{[r]},1^{[s]},1^{[t]})$ ($rs = 0$; $s < t$)  } \\ \cline{1-3}

\multicolumn{3}{l}{In $M_1 = \bar{A}_1 A_1$ $(p \geq 5)$ $(C_3(\#\cthree{7}))$} \\
\cthree{12} & \multicolumn{2}{l}{ $A_1$  via $(1^{[r]},1^{[s]})$ $(rs=0)$ } \\
\hline 
\multicolumn{3}{l}{In $M_1 = \bar{A}_1^3$ $(C_3(\#\cthree{8}))$} \\
\cthree{13} & \multicolumn{2}{l}{ $\bar{A}_1 A_1$  via $(1_a,1_a^{[r]},1_b)$ $(r \neq 0)$ } \\
\cthree{14} & \multicolumn{2}{l}{ $A_1$ via $(1,1^{[r]},1^{[s]})$  $(0 < r < s)$ } \\ \hline \hline
\multicolumn{3}{l}{In $M = A_1 A_1$ $(p \neq 2)$ $(C_3(\#\cthree{2}))$} \\
\cthree{15} & \multicolumn{2}{l}{ $A_1$  via $(1^{[r]},1^{[s]})$ $(rs=0; r \neq s)$ } \\
$\mathit{\cthree{12}^{\{\underline{0}\}}}$ & \multicolumn{2}{l}{ $A_1$  via $(1,1)$ $(p \geq 5)$} \\ 
\hline \hline
\multicolumn{3}{l}{In $M = G_2$ $(p = 2)$ $(C_3(\#\cthree{4}))$} \\ 
$\mathit{\cthree{9}^{\{\delta_1\}}}$ & $\tilde{A}_1 A_1$ & \\ 
\hline \hline
\end{longtable}

\begin{longtable}[c]{lll}
\caption{Irreducible subgroups of $L_{0} = B_3$. \label{B3tab}} \\
\hline
ID & Subgroup $X$ & $V_{L_{0}} \downarrow X$ \\
\hline

\bthree{1}a & $\bar{A}_1^2 \tilde{A}_1$ ($p \neq 2$) & $(1,1,0) + (0,0,2)$  \\

\bthree{2} & $\bar{A}_3$ & $010 + 0$ $(p \neq 2)$ or $010$ $(p=2)$ \\

\bthree{3} & $G_2$ & $10$  \\

\bthree{4} & $B_2 \tilde{A}_1$ $(p=2)$ & $(10,0) + (0,2)$ \\

\hline \hline
\endfirsthead
\multicolumn{3}{l}{In $M = \bar{A}_1^2 \tilde{A}_1$ $(B_3(\#\bthree{1}))$ } \\

\bthree{5} & \multicolumn{2}{l}{ $A_1 \tilde{A}_1$ via $(1_a,1_a^{[r]},1_b)$ (if $p=2$ then $r \neq 0)$ } \\

\bthree{6} & \multicolumn{2}{l}{ $\bar{A}_1 A_1$ via $(1_a^{[r]},1_b,1_a^{[s]})$ ($rs=0$) } \\

\bthree{7} & \multicolumn{2}{l}{$A_1$ via  $(1^{[r]},1^{[s]},1^{[t]})$ ($rt = 0$; $r \le s$; if $r=s$ then $p \neq 2$ and $r < t$)} \\

\hline \hline

\multicolumn{3}{l}{In $M = \bar{A}_3$ $(B_3(\#\bthree{2}))$} \\

$\mathit{\bthree{5}^{\{0\}}}$ & $\tilde{A}_1^2$ $(p \neq 2)$ & \\ 
\hline \hline

\multicolumn{3}{l}{In $M = G_2$ $(B_3(\#\bthree{3}))$} \\

$\mathit{\bthree{6}^{\{\underline{0}\}}}$ & $\bar{A}_1 A_1$ & \\
\bthree{8} & $A_2$ $(p=3)$ & $11$ \\

\bthree{9} & $A_1$ $(p \geq 7)$ & $6$ \\

\hline

\multicolumn{3}{l}{In $M_1 = A_2$ $(p=3)$ $(B_3(\#\bthree{8}))$} \\

$\mathit{\bthree{7}^{\{\delta_2\}}}$ & $A_1$ &  \\

\hline \hline

\multicolumn{3}{l}{In $M = B_2 \tilde{A}_1$ $(p=2)$ $(B_3(\#\bthree{4}))$} \\

\bthree{1}b & $\bar{A}_1^2 \tilde{A}_1$ & $(1,1,0) + (0,0,2)$ \\

\bthree{10} & $\tilde{A}_1^3$  & $(2,0,0) + (0,2,0) + (0,0,2)$ \\

\hline

\multicolumn{3}{l}{In $M_1 = \tilde{A}_1^3$ $(p=2)$ $(B_3(\#\bthree{10}))$} \\

\bthree{11} & \multicolumn{2}{l}{ $\tilde{A}_1 A_1$  via $(1_a,1_a^{[r]},1_b)$ $(r \neq 0)$  } \\

\bthree{12} & \multicolumn{2}{l}{ $A_1$  via $(1,1^{[r]},1^{[s]})$ $(0 < r < s)$ } \\

\hline \hline
\end{longtable}

\begin{longtable}{lll}
\caption{Irreducible subgroups of $L_{0} = D_4$. \label{D4tab}} \\
\hline
ID & Subgroup $X$ & $V_{L_{0}} \downarrow X$ \\
\hline

\dfour{1} & $\bar{A}_1^4$ & $(1,1,0,0) + (0,0,1,1)$ \\

\dfour{2} & $B_3$ & $T(\lambda_1;0)$ \\

\dfour{3} & $B_3$ & $001$ \\

\dfour{4} & $B_3$ & $001$ \\
\endfirsthead
\dfour{5}a & $A_1 B_2$ $(p \neq 2)$ & $(2,0) + (0,10)$ \\

\dfour{6}a & $A_1 B_2$ $(p \neq 2)$ & $(1,01)$ \\

\dfour{7}a & $A_1 B_2$ $(p \neq 2)$ & $(1,01)$ \\

\dfour{8}a & $A_2$ $(p \neq 3)$ & $11$ \\

\hline \hline

\multicolumn{3}{l}{In $M = \bar{A}_1^4$ ($D_4(\#\dfour{1})$)} \\

\dfour{9} & \multicolumn{2}{l}{ $A_1 \bar{A}_1^2$ via $(1_a,1_a^{[r]},1_b,1_c)$   } \\

\dfour{10} & \multicolumn{2}{l}{ $A_1 \bar{A}_1^2$ via $(1_a,1_b,1_a^{[r]},1_c)$   } \\

\dfour{11} & \multicolumn{2}{l}{ $A_1 \bar{A}_1^2$ via $(1_a,1_b,1_c,1_a^{[r]})$   } \\

\dfour{12} & \multicolumn{2}{l}{ $A_1 \bar{A}_1$ via $(1_a^{[r]},1_a^{[s]},1_a^{[t]},1_b)$ $(rst=0)$ } \\

\dfour{13} & \multicolumn{2}{l}{ $A_1 A_1$ via $(1_a,1_a^{[r]},1_b^{[s]},1_b^{[t]})$ $(st=0; r \leq s+t;$ if $r=0$ then $t \neq 0$)  } \\

\dfour{14} & \multicolumn{2}{l}{ $A_1 A_1$ via $(1_a,1_b^{[s]},1_a^{[r]},1_b^{[t]})$ $(st=0; r \leq s+t;$ if $r=0$ then $t \neq 0$) } \\

\dfour{15} & \multicolumn{2}{l}{ $A_1 A_1$ via $(1_a,1_b^{[s]},1_b^{[t]},1_a^{[r]})$ $(st=0; r \leq s+t;$ if $r=0$ then $t \neq 0$)  } \\

\dfour{16} & \multicolumn{2}{l}{ $A_1$ via $(1,1^{[r]},1^{[s]},1^{[t]})$ (if $r=0$ then $s < t$; if $s=0$ then $r < t$; if $t=0$ then $r < s$)  } \\

\hline \hline

\multicolumn{3}{l}{In $M = B_3$ ($D_4(\#\dfour{2})$) } \\

$\mathit{\dfour{9}^{\{0\}}}$ & $\bar{A}_1^2 A_1$ $(p \neq 2)$ & \\

\dfour{17} & $G_2$ & $T(10;0)$ \\ 
\dfour{5}b & $A_1 B_2$ $(p=2)$ & $\myTheta(2;10)$ \\ 
\hline

\multicolumn{3}{l}{In $M_1 = G_2$ ($D_4(\#\dfour{17})$)} \\

$\mathit{\dfour{12}^{\{\underline{0}\}}}$ & $\bar{A}_1 A_1$ & \\

\dfour{8}b & $A_2$ $(p=3)$ & $11 + 0$ \\

\dfour{18} & $A_1$ & $6 + 0$ \\

\hline \hline

\multicolumn{3}{l}{In $M = B_3$ ($D_4(\#\dfour{3})$)} \\

$\mathit{\dfour{11}^{\{0\}}}$ & $\bar{A}_1^2 A_1$ $(p \neq 2)$ &  \\

\it{\dfour{17}} & $G_2$ &  \\

\dfour{6}b & $A_1 B_2$ $(p=2)$ & $(1,01)$  \\

\hline \hline

\multicolumn{3}{l}{In $M = B_3$ ($D_4(\#\dfour{4})$)} \\

$\mathit{\dfour{10}^{\{0\}}}$ & $\bar{A}_1^2 A_1$ $(p \neq 2)$ &   \\

\it{\dfour{17}} & $G_2$ & \\

\dfour{7}b & $A_1 B_2$ $(p=2)$ & $(1,01)$ \\

\hline \hline

\multicolumn{3}{l}{In $M = A_1 B_2$ ($D_4(\#\dfour{5})$)} \\

$\mathit{\dfour{9}^{\{0\}}}$ & $A_1 \bar{A}_1^2$ &  \\

\dfour{19} & $A_1 A_1$ $( p \geq 5)$ & $(2,0) + (0,4)$ \\

\dfour{20} & $A_1^3$ $(p = 2)$ & $\myTheta(2;2;2)$  \\ 
\hline

\multicolumn{3}{l}{In $M_1 = A_1 A_1$ $(p \geq 5)$ ($D_4(\#\dfour{19})$)} \\

\dfour{21} & \multicolumn{2}{l}{ $A_1$ via $(1^{[r]},1^{[s]})$ $(rs=0)$ } \\

\hline \hline

\multicolumn{3}{l}{In $M = A_1^3$ $(p = 2)$ ($D_4(\#\dfour{20})$)} \\

\dfour{22} & \multicolumn{2}{l}{ $A_1 A_1$ via $(1_a,1_a^{[r]},1_b)$ $(r \neq 0)$ } \\

\dfour{23} & \multicolumn{2}{l}{ $A_1$ via $(1,1^{[r]},1^{[s]})$ $(0 < r < s)$ } \\

\hline \hline

\multicolumn{3}{l}{In $M = A_1 B_2$ ($D_4(\#\dfour{6})$)} \\

$\mathit{\dfour{11}^{\{0\}}}$ & $A_1 \bar{A}_1^2$ & \\ 
\dfour{24} & $A_1 A_1$ $(p \geq 5)$ & $(1,3)$ \\

\dfour{25} & $A_1^3$ $(p = 2)$ & $(1,1,1)$  \\

\hline

\multicolumn{3}{l}{In $M_1 = A_1 A_1$ $(p \geq 5)$ ($D_4(\#\dfour{24})$)} \\

\dfour{26} & \multicolumn{2}{l}{ $A_1$ via $(1^{[r]},1^{[s]})$ $(rs=0)$ } \\

\hline

\multicolumn{3}{l}{In $M_1 = A_1^3$ $(p = 2)$ ($D_4(\#\dfour{25})$)} \\

\dfour{27} & \multicolumn{2}{l}{ $A_1 A_1$ via $(1_a,1_a^{[r]},1_b)$ $(r \neq 0)$ } \\

\dfour{28} & \multicolumn{2}{l}{ $A_1$ via $(1,1^{[r]},1^{[s]})$ $(0 < r < s)$ } \\

\hline \hline

\multicolumn{3}{l}{In $M = A_1 B_2$ ($D_4(\#\dfour{7})$)} \\

$\mathit{\dfour{10}^{\{0\}}}$ & $A_1 \bar{A}_1^2$ & \\ 
\dfour{29} & $A_1 A_1$ $(p \geq 5)$ & $(1,3)$ \\

\dfour{30} & $A_1^3$ $(p = 2)$ & $(1,1,1)$  \\

\hline

\multicolumn{3}{l}{In $M_1 = A_1 A_1$ $(p \geq 5)$ ($D_4(\#\dfour{29})$)} \\

\dfour{31} & \multicolumn{2}{l}{ $A_1$ via $(1^{[r]},1^{[s]})$ $(rs=0)$ } \\

\hline

\multicolumn{3}{l}{In $M_1 = A_1^3$ $(p = 2)$ ($D_4(\#\dfour{30})$)} \\

\dfour{32} & \multicolumn{2}{l}{ $A_1 A_1$ via $(1_a,1_a^{[r]},1_b)$ $(r \neq 0)$ } \\

\dfour{33} & \multicolumn{2}{l}{ $A_1$ via $(1,1^{[r]},1^{[s]})$ $(0 < r < s)$ } \\

\hline \hline

\multicolumn{3}{l}{In $M = A_2$ ($D_4(\#\dfour{8})$)} \\

$\mathit{\dfour{21}^{\{\underline{0}\}}}$ & $A_1$ $(p \geq 5)$ & \\ 
$\mathit{\dfour{16}^{\{\delta_3\}}}$ & $A_1$ $(p=3)$ & \\ 
\hline \hline

\end{longtable}

\begin{longtable}{lll}
\caption{Irreducible subgroups of $L_{0} = D_5$. \label{D5tab}} \\
\hline
ID & Subgroup $X$ & $V_{L_{0}} \downarrow X$ \\
\hline

\dfive{1} & $\bar{A}_1^2 \bar{A}_3$ & $(1,1,0) + (0,0,010)$ \\

\dfive{2}a & $B_2^2$ $(p \neq 2)$ & $(10,0) + (0,10)$ \\

\dfive{3}a & $A_1 B_3$ $(p \neq 2)$ & $(2,0) + (0,100)$ \\

\dfive{4} & $B_4$ & $T(\lambda_1;0)$ \\ 
\dfive{5} & $B_2$ $(p \neq 2)$ & $02$ \\

\dfive{6} & $B_2$ $(p \neq 2)$ & $02$ \\

\hline \hline \endfirsthead

\multicolumn{3}{l}{In $M = \bar{A}_1^2 \bar{A}_3$ $(D_5(\#\dfive{1}))$} \\

\dfive{7} & $\bar{A}_1^2 B_2$ &  $(1,1,0) + (0,0,T(10;0))$ \\ 
\dfive{8}a & $\bar{A}_1^2 A_1^2$ $(p \neq 2)$ & $(1,1,0,0) + (0,0,2,0) +  (0,0,0,2)$ \\ 

\dfive{9} & \multicolumn{2}{l}{ $A_1 \bar{A}_3$ via $(1,1^{[r]},100)$ } \\

\hline

\multicolumn{3}{l}{In $M_1 = \bar{A}_1^2 B_2$ $(D_5(\#\dfive{7}))$} \\

\dfive{10} & $\bar{A}_1^2 A_1$ $(p \geq 5)$ & $(1,1,0) + (0,0,4) + 0$ \\

\dfive{8}b & $\bar{A}_1^2 A_1^2$ $(p=2)$ & $(1,1,0,0) + (0,0,\myTheta(2;2))$ \\ 
\dfive{11} & \multicolumn{2}{l}{ $A_1 B_2$ via $(1_a,1_a^{[r]},10)$ $(r \neq 0)$} \\

\hline

\multicolumn{3}{l}{In $M_2 = \bar{A}_1^2 A_1$ $(p \geq 5)$ $(D_5(\#\dfive{10}))$} \\

\dfive{12} & \multicolumn{2}{l}{ $A_1 A_1$ via $(1_a,1_a^{[r]},1_b)$ $(r \neq 0)$  } \\

\dfive{13} & \multicolumn{2}{l}{ $\bar{A}_1 A_1$ via $(1_a^{[r]},1_b,1_a^{[s]})$ $(rs=0)$ } \\

\dfive{14} & \multicolumn{2}{l}{ $A_1$ via $(1^{[r]},1^{[s]},1^{[t]})$ $(rt=0; r < s)$ } \\

\hline

\multicolumn{3}{l}{In $M_2 = \bar{A}_1^2 A_1^2$ $(D_5(\#\dfive{8}))$} \\
\dfive{15} & \multicolumn{2}{l}{ $A_1 A_1^2$ via $(1_a,1_a^{[r]},1_b,1_c)$ (if $p=2$ then $r \neq 0$) } \\

\dfive{16} & \multicolumn{2}{l}{ $A_1 \bar{A}_1 A_1$ via $(1_a^{[r]},1_b,1_a^{[s]},1_c)$ $(rs=0)$ } \\

\dfive{17} & \multicolumn{2}{l}{ $\bar{A}_1^2 A_1$ via $(1_a,1_b,1_c,1_c^{[r]})$ $(r \neq 0)$ } \\
\dfive{18} & \multicolumn{2}{l}{ $A_1 A_1$ via $(1_a^{[r]},1_a^{[s]},1_a^{[t]},1_b)$ ($rt=0$; $r \leq s$; if $r=s$ then $p \neq 2$ and $r < t$) } \\

\dfive{19} & \multicolumn{2}{l}{ $A_1 \bar{A}_1$ via $(1_a^{[r]},1_b,1_a^{[s]},1_a^{[t]})$ ($rs=0$; $s < t$) } \\
\dfive{20} & \multicolumn{2}{l}{ $A_1 A_1$ via $(1_a,1_a^{[r]},1_b,1_b^{[s]})$ (if $p=2$ then $r \neq 0$; $s \neq 0$) } \\

\dfive{21} & \multicolumn{2}{l}{ $A_1 A_1$ via $(1_a^{[r]},1_b^{[t]},1_a^{[s]},1_b^{[u]})$ ($rs=tu=0$; $r \leq t$; if $r=t$ then $s \leq u$) } \\
\dfive{22} & \multicolumn{2}{l}{ $A_1$ via $(1^{[r]},1^{[s]},1^{[t]},1^{[u]})$ ($rt=0$; $r \leq s$; $t < u$; if $r=s$ then $p \neq 2$ and $r < t$) } \\ 
\hline \hline 

\multicolumn{3}{l}{In $M = B_2^2$ $(D_5(\#\dfive{2}))$ } \\

\it{\dfive{7}} & $\bar{A}_1^2 B_2$  & \\

\dfive{23} & $A_1 B_2$ $(p \geq 5)$ & $(4,0) + (0,10)$ \\

\dfive{24} & $A_1^2 B_2$ $(p=2)$ & $\myTheta(2;2;10)$ \\ 
\dfive{25} & \multicolumn{2}{l}{ $B_2$ via $(10,10^{[r]})$ $(r \neq 0)$ } \\

\dfive{26} & \multicolumn{2}{l}{ $B_2$ via $(10,02^{[r]})$ $(p=2)$ } \\

\hline

\multicolumn{3}{l}{In $M_1 = A_1 B_2$ $(p \geq 5)$ $(D_5(\#\dfive{23}))$} \\

\it{\dfive{10}} & $A_1 \bar{A}_1^2$ & \\
\dfive{27} & $A_1^2$ & $(4,0) + (0,4)$ \\

\hline

\multicolumn{3}{l}{In $M_2 = A_1^2$ $(p \geq 5)$ $(D_5(\#\dfive{27}))$} \\

\dfive{28} & \multicolumn{2}{l}{ $A_1$ via $(1,1^{[r]})$ $(r \neq 0)$ } \\

\hline

\multicolumn{3}{l}{In $M_1 = A_1^2 B_2$ $(p=2)$ $(D_5(\#\dfive{24}))$} \\

\it{\dfive{8}} & $A_1^2 \bar{A}_1^2$ & \\

\dfive{29} & $A_1^4$ & $\myTheta(2;2;2;2)$ \\ 
\dfive{30} & \multicolumn{2}{l}{ $A_1 B_2$ via $(1_a,1_a^{[r]},10)$ $(r \neq 0)$} \\

\hline

\multicolumn{3}{l}{In $M_2 = A_1^4$ $(p=2)$ $(D_5(\#\dfive{29}))$ } \\

\dfive{31} & \multicolumn{2}{l}{ $A_1 A_1^2$ via $(1_a,1_a^{[r]},1_b,1_c)$ $(r \neq 0)$ } \\

\dfive{32} & \multicolumn{2}{l}{ $A_1 A_1$ via $(1_a,1_a^{[r]},1_a^{[s]},1_b)$ $(0 < r < s)$ } \\

\dfive{33} & \multicolumn{2}{l}{ $A_1 A_1$ via $(1_a,1_a^{[r]},1_b,1_b^{[s]})$ $(rs \neq 0; r \leq s)$ } \\

\dfive{34} & \multicolumn{2}{l}{ $A_1$ via $(1,1^{[r]},1^{[s]},1^{[t]})$ $(0 < r < s < t)$ } \\

\hline \hline

\multicolumn{3}{l}{In $M = A_1 B_3$ $(D_5(\#\dfive{3}))$ } \\

$\mathit{\dfive{9}^{\{0\}}}$ & $A_1 \bar{A}_3$ &  \\

\it{\dfive{8}} & $A_1^2 \bar{A}_1^2$ $(p \neq 2)$ &  \\

\dfive{35} & $A_1 G_2$ & $(2,0) + (0,10)$ $(p \neq 2)$ or $\myTheta(2;10)$ $(p=2)$ \\ 
\it{\dfive{24}} & $A_1^2 B_2$ $(p=2)$ &  \\

\hline

\multicolumn{3}{l}{In $M_1 = A_1 G_2$ $(D_5(\#\dfive{35}))$ } \\

$\mathit{\dfive{16}^{\{\underline{0}\}}}$ & $A_1 \bar{A}_1 A_1$ &  \\

\dfive{36} & $A_1 A_2$ $(p=3)$ & $(2,0) + (0,11)$ \\

\dfive{37} & $A_1 A_1$ $(p \geq 7)$ & $(2,0) + (0,6)$  \\

\hline

\multicolumn{3}{l}{In $M_2 = A_1 A_2$ $(p=3)$ $(D_5(\#\dfive{36}))$ } \\

$\mathit{\dfive{18}^{\{\delta_2\}}}$ & $A_1 A_1$ & \\

\hline

\multicolumn{3}{l}{In $M_2 = A_1 A_1$ $(p \geq 7)$ $(D_5(\#\dfive{37}))$ } \\

\dfive{38} & \multicolumn{2}{l}{ $A_1$ via $(1^{[r]},1^{[s]})$ $(rs=0)$ } \\

\hline \hline

\multicolumn{3}{l}{In $M = B_4$ $(D_5(\#\dfive{4}))$} \\

$\mathit{\dfive{9}^{\{0\}}}$ & $A_1 \bar{A}_3$ $(p \neq 2)$ &  \\

\it{\dfive{7}} & $\bar{A}_1^2 B_2$ $(p \neq 2)$ &  \\

\dfive{39} & $A_1^2$ $(p \neq 2)$ & $(2,2) + 0$  \\

\dfive{40} & $A_1$ $(p \geq 11)$ & $8 + 0$  \\

\dfive{2}b & $B_2^2$ $(p=2)$ & $\myTheta(10;10)$  \\ 
\dfive{3}b & $A_1 B_3$ $(p=2)$ & $\myTheta(2;100)$  \\ 
\hline

\multicolumn{3}{l}{In $M_1 = A_1^2$ $(p \neq 2)$ $(D_5(\#\dfive{39}))$} \\

\dfive{41}& \multicolumn{2}{l}{ $A_1$ via $(1,1^{[r]})$ ($r \neq 0$) } \\

\hline \hline

\multicolumn{3}{l}{In $M = B_2$ $(p \neq 2)$  $(D_5(\#\dfive{5}))$} \\

$\mathit{\dfive{21}^{\{\underline{0}\}}}$ & $A_1^2$ & \\

$\mathit{\dfive{38}^{\{0\}}}$ & $A_1$ $(p \geq 7)$ & \\

\hline \hline

\multicolumn{3}{l}{In $M = B_2$ $(p \neq 2)$  $(D_5(\#\dfive{6}))$} \\

$\mathit{\dfive{21}^{\{\underline{0}\}}}$ & $A_1^2$ & \\

$\mathit{\dfive{38}^{\{0\}}}$ & $A_1$ $(p \geq 7)$ & \\

\hline \hline

\end{longtable}

} 
In the following table, subgroups of $D_6$ which are not centralised by a graph automorphism can be distinguished by their actions on the two half-spin modules $V_{D_6}(\lambda_5)$ and $V_{D_6}(\lambda_6)$. For consistency, in the following table we fix how the subgroups $A_1 C_3 = D_6(\#6)$, $A_1 C_3 = D_6(\#7)$ and $\bar{A}_1^2 A_1 B_2 = D_6(\#12)$ act on these modules.  

We take $A_1C_3 = D_6(\#\dsix{6})$ to act on $V_{D_6}(\lambda_5)$ with factors $(W(3),0)/(1,W(010))$ and therefore act on $V_{D_6}(\lambda_6)$ with factors $(W(2),100)/(0,W(001))$. Taking an image under a graph automorphism, it follows that $A_1C_3 = D_6(\#\dsix{7})$ acts on $V_{D_6}(\lambda_5)$ with factors $(W(2),100)/(0,W(001))$ and on $V_{D_6}(\lambda_6)$ with factors $(W(3),0)/(1,W(010))$. 

We also take $\bar{A}_1^2 A_1 B_2 = D_6(\#\dsix{12})$ to act on $V_{D_6}(\lambda_5)$ with factors $(1,0,W(2),0)/\allowbreak(1,0,0,W(10))/\allowbreak(0,1,1,01)$ and on $V_{D_6}(\lambda_6)$ with factors $(1,0,1,01)/\allowbreak(0,1,W(2),0)/\allowbreak(0,1,0,W(10))$. This also determines the action of the subgroups $\bar{A}_1^2 A_1 A_1 = D_6(\#\dsix{73})$ $(p \geq 5)$ and $\bar{A}_1^2 A_1 A_1 = D_6(\#\dsix{74})$ $(p \geq 5)$.

{
\small 
\begin{longtable}[c]{lll}
\caption{Irreducible subgroups of $L_{0} = D_6$. \label{D6tab}} \\
\hline
ID & Subgroup $X$ & $V_{L_{0}} \downarrow X$ \\
\hline

\dsix{1} & $\bar{A}_1^2 D_4$ & $(1,1,0) + (0,0,\lambda_1)$ \\

\dsix{2} & $\bar{A}_3^2$ & $(010,0) + (0,010)$ \\

\dsix{3}a & $A_1 B_4$ $(p \neq 2)$ & $(2,0) + (0,\lambda_1)$ \\

\dsix{4}a & $B_2 B_3$ $(p \neq 2)$ & $(10,0) + (0,100)$ \\

\dsix{5} & $B_5$ & $T(\lambda_1;0)$ \\ 
\dsix{6} & $A_1 C_3$ & $(1,100)$ \\
\hline \hline
\end{longtable}
}
We remind the reader that in the following table, subgroups are listed up to $\Aut(L_{0})$-conjugacy. 
{\small

\begin{longtable}[c]{ll>{\raggedright\arraybackslash}p{0.6\textwidth-\tabcolsep}@{}}
\caption{Irreducible subgroups of $L_{0} = D_7$. \label{D7tab}} \\
\hline
ID & Subgroup $X$ & $V_{L_{0}} \downarrow X$ \\
\hline

\dseven{1} & $\bar{A}_1^2 D_5$ & $(1,1,0) + (0,0,\lambda_1)$ \\

\dseven{2} & $\bar{A}_3 D_4$ & $(010,0) + (0,\lambda_1)$ \\

\dseven{3} & $B_6$ & $T(\lambda_1;0)$ \\ 

\dseven{4}a & $B_3^2$ $(p \neq 2)$ & $(100,0) + (0,100)$ \\

\dseven{5}a & $B_2 B_4$ $(p \neq 2)$ & $(10,0) + (0,\lambda_1)$ \\

\dseven{6}a & $A_1 B_5$ $(p \neq 2)$ & $(2,0) + (0,\lambda_1)$ \\
\hline \hline

\end{longtable}

}
\vspace{-0.39cm}
\subsection{Irreducible diagonal subgroups} \label{tabs:diags}

In this section we give the tables of diagonal subgroups referred to in Tables \ref{D6tab} and \ref{D7tab}. The first column gives the ID number, as in the previous tables, and the second column gives the embedding of the diagonal subgroups. To describe the embeddings we use a slightly modified notation, to shorten the tables. Specifically, we introduce a shorthand for diagonal subgroups of $A_1^n Z$, where $Z$ has no simple factor of type $A_1$. For example, instead of writing $A_1^2 B_2 \hookrightarrow A_1^3 B_2$ via $(1_a^{[r]},1_a^{[s]},1_b,10)$ we just write $(a^{[r]},a^{[s]},b,10)$; from any such vector it is easy to recover the isomorphism type of the diagonal subgroup. Moreover, in Tables \ref{d7diag13} and \ref{d7diag279} the usual notation for diagonal subgroups is used but we again omit the isomorphism type of each diagonal subgroup as they too can be easily recovered from the listed embedding. 

Table \ref{condD639} gives the extra restrictions on the field twists for $D_6(\#\dsix{38})$, a diagonal subgroup of $\bar{A}_1^6 = D_6(\#\dsix{8})$. The restrictions are given in rows of the tables: the first column lists all permitted equalities amongst the field twists; the second column lists any further requirements. So an ordered set $\{r, \ldots , w\}$ is permitted if it satisfies the conditions in the first and second column of a row of the table. We note that a set of field twists satisfies the conditions of at most one row. We emphasise that an ordered set may be excluded either because it yields a $D_7$-reducible subgroup, or because it yields a repeated diagonal subgroup.  

\subsubsection{Diagonal subgroups contained in $D_6$} \leavevmode

{ \small 
\begin{longtable}{>{\raggedright\arraybackslash}p{0.05\textwidth - 2\tabcolsep}>{\raggedright\arraybackslash}p{0.95\textwidth-\tabcolsep}@{}}
\caption{Irreducible diagonal subgroups of $\bar{A}_1^6 = D_6(\#\dsix{8})$. \label{d6diag8}} \\
\hline
ID & Diagonal embedding \\
\hline
\dsix{15} & $(a,a^{[r]},b,c,d,e)$  \\
\dsix{16} & $(a,b,a^{[r]},c,d,e)$  \\
\dsix{17} &  $(a^{[r]},a^{[s]},a^{[t]},b,c,d)$ $(rt=0; r \leq s)$ \\
\dsix{18} & $(a,b,a^{[r]},c,a^{[s]},d)$ $(r \leq s)$ \\
\dsix{19} &  $(a,b,a^{[r]},c,d,a^{[s]})$ $(r \leq s)$ \\
\endfirsthead
\dsix{20} & $(a,a^{[r]},b,b^{[s]},c,d)$ $(r \leq s$; if $r=0$ then $s \neq 0$) \\
\dsix{21} & $(a^{[r]},a^{[s]},b,c,b^{[t]},d)$ $(rs=0)$ \\
\dsix{22} & $(a,b^{[s]},a^{[r]},b^{[t]},c,d)$ $(st=0; r \leq s+t$; if $r=0$ then $t \neq 0$) \\
\dsix{23} & $(a^{[r]},b^{[t]},a^{[s]},c,b^{[u]},d)$ $(rs=tu=0)$ \\
\dsix{24} & $(a,a^{[r]},a^{[s]},a^{[t]},b,c)$ $(s \leq t$; if $r=0$ then $s < t$; if $s=0$ then $r < t$; if $t=0$ then $r < s$) \\
\dsix{25} & $(a^{[r]},a^{[s]},a^{[t]},b,a^{[u]},c)$ $(rst=0; t \leq u);$ \\
\dsix{26} & $(a^{[r]},a^{[s]},a^{[t]},b^{[u]},b^{[v]},c)$ $(rst=uv=0)$ \\
\dsix{27} & $(a^{[r]},a^{[s]},a^{[t]},b,c,c^{[u]})$ ($rst=0$; if $u=0$ then $r < s$) \\
\hline
\end{longtable}

\begin{longtable}[c]{ll}
\caption{Conditions on field twists for $D_6(\#\dsix{38})$. \label{condD639}} \\

\hline \noalign{\smallskip}

Equalities among $r, \ldots, w$ & Further requirements on $r, \ldots, w$ \\

\hline 

none & $r < \text{min}\{s,t\}$; $t < \text{min}\{u,v,w\}$; $s < u$ \\
$r=s$ & $t < \text{min}\{u,v\}$; $v < w$ \\
$r=t$ & $s < u$ \\
\endfirsthead
$r=s=t$ & $v < w$ \\
$r=t=v$ & $s < u < w$ \\
$r=t=w$ & $s < u < v$ \\
$r=s$; $t=v$ & $u < w$ \\
$r=t$; $s=v$ & none \\
$r=s=t=v$ & $u < w$ \\
$r=s=t$; $u=v$ & none \\
$r=t$; $s=v$; $u=w$ & none \\
\hline
\end{longtable}

\begin{longtable}[c]{llll}
\caption{Irreducible diagonal subgroups of $\bar{A}_1^2 A_1 A_1 = D_6(\#\dsix{43})$ $(p \geq 5)$. \label{d6diag44}} \\
\hline
ID & Diagonal embedding & ID & Diagonal embedding \\
\hline
\dsix{48} & $(a,a^{[r]},b,c)$ & \dsix{53} & $(a^{[r]},a^{[s]},b,a^{[t]})$ ($rt=0; r \leq s$) \\
\dsix{49} & $(a^{[r]},b,a^{[s]},c)$ ($rs=0$) & \dsix{54} & $(a^{[r]},b,a^{[s]},a^{[t]})$ ($rst=0$) \\
\dsix{50} & $(a^{[r]},b,c,a^{[s]})$ ($rs=0$) & \dsix{55} & $(a,a^{[r]},b^{[s]},b^{[t]})$ $(st=0)$ \\
\dsix{51} & $(a,b,c^{[r]},c^{[s]})$ ($rs=0$) & \dsix{56} & $(a^{[r]},b^{[t]},a^{[s]},b^{[u]})$ ($rs=tu=0$) \\
\dsix{52} & $(a^{[r]},a^{[s]},a^{[t]},b)$ $(rt=0; r \leq s$; if $r=s$ then $r < t$) & \dsix{57} & $(a^{[r]},a^{[s]},a^{[t]},a^{[u]})$ $(rtu=0; r \leq s$; \\
& & & if $r=s$ then $r < t$) \\
\hline
\end{longtable}

\begin{longtable}[c]{ll}
\caption{Irreducible diagonal subgroups of $\bar{A}_1^2 A_1^3 = D_6(\#\dsix{44})$ $(p=2)$. \label{d6diag45}} \\
\hline
ID & Diagonal embedding \\
\hline
\dsix{58} & $(a,a^{[r]},b,c,d)$ $(r \neq 0)$ \\
\dsix{59} & $(a^{[r]},b,a^{[s]},c,d)$ $(rs=0)$ \\
\dsix{60} & $(a,b,c,c^{[r]},d)$ $(r \neq 0)$ \\
\dsix{61} & $(a^{[r]},a^{[s]},a^{[t]},b,c)$ $(rt=0; r < s)$ \\
\dsix{62} & $(a^{[r]},b,a^{[s]},a^{[t]},c)$ $(rs=0; s < t)$ \\
\dsix{63} & $(a,b,c,c^{[r]},c^{[s]})$ $(0 < r < s)$ \\
\dsix{64} & $(a,a^{[r]},b,b^{[s]},c)$ $(rs \neq 0)$ \\
\dsix{65} & $(a^{[r]},b^{[t]},a^{[s]},b^{[u]},c)$ $(rs=tu=0; r \leq t;$ if $r=t$ then $s \leq u$) \\
\dsix{66} & $(a^{[r]},b,a^{[s]},c,c^{[t]})$ $(rs=0 \neq t)$ \\
\dsix{67} & $(a^{[r]},a^{[s]},a^{[t]},a^{[u]},b)$ $(rt=0; r < s; t < u)$ \\
\dsix{68} & $(a^{[r]},b,a^{[s]},a^{[t]},a^{[u]})$ $(rs=0; s < t < u)$ \\
\dsix{69} & $(a^{[r]},a^{[s]},a^{[t]},b,b^{[u]})$ $(rt=0 \neq u; r < s)$ \\
\dsix{70} & $(a^{[r]},b^{[u]},a^{[s]},a^{[t]},b^{[v]})$ $(rs=uv=0; s < t)$ \\
\dsix{71} & $(a,a^{[r]},b,b^{[s]},b^{[t]})$ $(r \neq 0; 0 < s < t)$ \\
\dsix{72} & $(a^{[r]},a^{[s]},a^{[t]},a^{[u]},a^{[v]})$ $(rt=0; r < s; t < u < v)$ \\
\hline
\end{longtable}

\begin{longtable}{llll}
\caption{Irreducible diagonal subgroups of $\bar{A}_1^2 A_1 A_1 = D_6(\#\dsix{73})$ $(p \geq 5)$. \label{d6diag74}} \\
\hline
ID & Diagonal embedding & ID & Diagonal embedding \\
\hline
\dsix{79} & $(a^{[r]},a^{[s]},b,c)$ $(rs=0)$ & \dsix{86} & $(a^{[r]},a^{[s]},b,a^{[t]})$ ($rst=0$) \\
\dsix{80} & $(a^{[r]},b,a^{[s]},c)$ ($rs=0$) & \dsix{87} & $(a^{[r]},b,a^{[s]},a^{[t]})$ ($rst=0; s \neq t)$ \\
\dsix{81} & $(a^{[r]},b,c,a^{[s]})$ ($rs=0$) & \dsix{88} & $(a,b^{[r]},b^{[s]},b^{[t]})$ ($rst=0; s \neq t$) \\
\dsix{82} & $(a,b^{[r]},b^{[s]},c)$ ($rs=0$) & \dsix{89} & $(a^{[r]},a^{[s]},b^{[t]},b^{[u]})$ $(rs=tu=0; t \neq u)$ \\
\dsix{83} & $(a,b^{[r]},c,b^{[s]})$ ($rs=0$) & \dsix{90} & $(a^{[r]},b^{[t]},a^{[s]},b^{[u]})$ ($rs=tu=0$) \\
\dsix{84} & $(a,b,c^{[r]},c^{[s]})$ ($rs=0; r \neq s$) & \dsix{91} & $(a^{[r]},b^{[t]},b^{[u]},a^{[s]})$ ($rs=tu=0$) \\
\dsix{85} & $(a^{[r]},a^{[s]},a^{[t]},b)$ $(rst=0)$ & \dsix{92} & $(a^{[r]},a^{[s]},a^{[t]},a^{[u]})$ $(rstu=0; t \neq u)$ \\
\hline
\end{longtable}

\begin{longtable}{llll}
\caption{Irreducible diagonal subgroups of $\bar{A}_1^2 A_1^3 = D_6(\#\dsix{74})$ $(p=2)$. \label{d6diag75}} \\
\hline
ID & Diagonal embedding & ID & Diagonal embedding \\
\hline
\dsix{93} & $(a^{[r]},a^{[s]},b,c,d)$ $(rs=0; r \neq s)$ & \dsix{103} & $(a^{[r]},b,a^{[s]},c,c^{[t]})$ $(rs=0 \neq t)$ \\
\dsix{94} & $(a^{[r]},b,a^{[s]},c,d)$ $(rs=0)$ & \dsix{104} & $(a,b^{[r]},b^{[s]},c,c^{[t]})$ $(rs=0 \neq t)$ \\
\dsix{95} & $(a,b^{[r]},b^{[s]},c,d)$ $(rs=0)$ & \dsix{105} & $(a^{[r]},a^{[s]},a^{[t]},a^{[u]},b)$ $(rst=0; r \neq s; t < u)$ \\
\endfirsthead
\dsix{96} & $(a,b,c,c^{[r]},d)$ $(r \neq 0)$ & \dsix{106} & $(a^{[r]},b,a^{[s]},a^{[t]},a^{[u]})$ $(rs=0; s < t < u)$ \\
\dsix{97} & $(a^{[r]},a^{[s]},a^{[t]},b,c)$ $(rst=0; r \neq s)$ & \dsix{107} & $(a,b^{[r]},b^{[s]},b^{[t]},b^{[u]})$ $(rs=0; s < t < u)$ \\
\dsix{98} & $(a^{[r]},b,a^{[s]},a^{[t]},c)$ $(rs=0; s < t)$ & \dsix{108} & $(a^{[r]},a^{[s]},a^{[t]},b,b^{[u]})$ $(rst=0; r \neq s; u \neq 0)$ \\
\dsix{99} & $(a,b^{[r]},b^{[s]},b^{[t]},c)$ $(rs=0; s < t)$ & \dsix{109} & $(a^{[r]},b^{[u]},a^{[s]},a^{[t]},b^{[v]})$ $(rs=uv=0; s < t)$ \\
\dsix{100} & $(a,b,c,c^{[r]},c^{[s]})$ $(0 < r < s)$ & \dsix{110} & $(a^{[r]},b^{[t]},b^{[u]},b^{[v]},a^{[s]})$ $(rs=tu=0; u < v)$ \\
\dsix{101} & $(a^{[r]},a^{[s]},b,b^{[t]},c)$ $(rs=0; r \neq s; t \neq 0)$ & \dsix{111} & $(a^{[r]},a^{[s]},b,b^{[t]},b^{[u]})$ $(rs=0; r \neq s; 0 < t < u)$ \\
\dsix{102} & $(a^{[r]},b^{[t]},a^{[s]},b^{[u]},c)$ $(rs=tu=0$) & \dsix{112} & $(a^{[r]},a^{[s]},a^{[t]},a^{[u]},a^{[v]})$ $(rst=0; r \neq s; t < u < v)$ \\
\hline
\end{longtable}

\begin{longtable}{llll}
\caption{Irreducible diagonal subgroups of $A_1^5 = D_6(\#\dsix{137})$ $(p = 2)$. \label{d6diag138}} \\
\hline
ID & Diagonal embedding & ID & Diagonal embedding\\
\hline
\dsix{140} & $(a,a^{[r]},b,c,d)$ $(r \neq 0)$ & \dsix{143} & $(a,a^{[r]},a^{[s]},a^{[t]},b)$ $(0 < r < s < t)$ \\
\dsix{141} & $(a,a^{[r]},a^{[s]},b,c)$ $(0 < r < s)$ & \dsix{144} & $(a,a^{[r]},a^{[s]},b,b^{[t]})$ $(0 < r < s; t \neq 0)$ \\
\dsix{142} & $(a,a^{[r]},b,b^{[s]},c)$ $(rs \neq 0; r \leq s)$ & \dsix{145} & $(a,a^{[r]},a^{[s]},a^{[t]},a^{[u]})$ $(0 < r < s < t < u)$ \\
\hline
\end{longtable}

\subsubsection{Diagonal subgroups contained in $D_7$} \leavevmode

\begin{longtable}{ll}
\caption{Irreducible diagonal subgroups of $\bar{A}_1^4 \bar{A}_3 = D_7(\#\dseven{10})$. \label{d7diag10}} \\
\hline
ID & Diagonal embedding \\
\hline

\dseven{18} & $(a,a^{[r]},b,c,100)$ \\

\dseven{19} & $(a,b,a^{[r]},c,100)$ \\

\dseven{20} & $(a^{[r]},a^{[s]},a^{[t]},b,100)$ $(rt=0; r \leq s)$ \\

\dseven{21} & $(a,a^{[r]},b,b^{[s]},100)$ $(r \leq s$; if $r=0$ then $s \neq 0$) \\ 
\dseven{22} & $(a,b^{[s]},a^{[r]},b^{[t]},100)$ $(st=0; r \leq s+t$; if $r=0$ then $t \neq 0$) \\

\dseven{23} & $(a,a^{[r]},a^{[s]},a^{[t]},100)$ $(s \leq t$; if $r=0$ then $s < t$; if $s=0$ then $r < t$) \\

\hline
\end{longtable}

\begin{longtable}{ll}
\caption{Irreducible diagonal subgroups of $\bar{A}_1^4 B_2 = D_7(\#\dseven{16})$. \label{d7diag16}} \\
\hline
ID & Diagonal embedding \\
\hline

\dseven{25} & $(a,a^{[r]},b,c,10)$ $(r \neq 0)$ \\

\dseven{26} & $(a,b,a^{[r]},c,10)$ \\

\dseven{27} & $(a^{[r]},a^{[s]},a^{[t]},b,10)$ $(rt=0; r < s)$ \\

\dseven{28} & $(a,a^{[r]},b,b^{[s]},10)$ $(rs \neq 0; r \leq s)$ \\

\dseven{29} & $(a,b^{[s]},a^{[r]},b^{[t]},10)$ $(st=0; r \leq s+t$; if $r=0$ then  $t \neq 0$) \\

\dseven{30} & $(a,a^{[r]},a^{[s]},a^{[t]},10)$ $(r \neq 0; s < t$; if $s=0$ then $r < t$) \\

\hline
\end{longtable}
\vspace{-0.1cm}
\begin{longtable}{ll}
\caption{Irreducible diagonal subgroups of $\bar{A}_1^4 A_1 = D_7(\#\dseven{24})$ $(p \geq 5)$. \label{d7diag24}} \\
\hline
ID & Diagonal embedding \\
\hline
\dseven{31} & $(a,a^{[r]},b,c,d)$ $(r \neq 0)$ \\

\dseven{32} & $(a,b,a^{[r]},c,d)$ \\

\dseven{33} & $(a^{[r]},b,c,d,a^{[s]})$ $(rs=0)$ \\
\endfirsthead
\dseven{34} & $(a^{[r]},a^{[s]},a^{[t]},b,c)$ $(rt=0; r < s)$ \\

\dseven{35} & $(a^{[r]},a^{[s]},c,d,a^{[t]})$ $(rt=0; r < s)$ \\

\dseven{36} & $(a^{[r]},b,a^{[s]},c,a^{[t]})$ $(rt=0; r \leq s)$ \\
\dseven{37} & $(a,a^{[r]},b,b^{[s]},c)$ $(rs \neq 0; r \leq s)$ \\

\dseven{38} & $(a,a^{[r]},b^{[s]},c,b^{[t]})$ $(r \neq 0; st=0)$ \\

\dseven{39} & $(a,b^{[s]},a^{[r]},b^{[t]},c)$ $(st=0; r \leq s+t$; if $r=0$ then $t \neq 0$) \\

\dseven{40} & $(a^{[r]},b^{[t]},a^{[s]},c,b^{[u]})$ $(rs=tu=0)$ \\
\dseven{41} & $(a,a^{[r]},a^{[s]},a^{[t]},b)$ $(r \neq 0; s < t$; if $s=0$ then $r < t$) \\

\dseven{42} & $(a^{[r]},a^{[s]},a^{[t]},b,a^{[u]})$ $(rtu=0; r < s)$ \\

\dseven{43} & $(a,a^{[r]},b^{[s]},b^{[t]},b^{[u]})$ $(r \neq 0; su=0; s < t)$ \\

\dseven{44} & $(a,b^{[s]},a^{[r]},b^{[t]},b^{[u]})$ $(stu=0$; if $r=0$ then $s \leq t$) \\

\dseven{45} & $(a^{[r]},b^{[t]},b^{[u]},b^{[v]},a^{[s]})$ $(rs=tu=0; u < v)$ \\

\dseven{46} & $(a^{[r]},a^{[s]},a^{[t]},a^{[u]},a^{[v]})$ ($rv=0; r < s; t < u; r \leq t$; if $r=t$ then $s < u$) \\

\hline
\end{longtable}
\vspace{-0.1cm}
\begin{longtable}{>{\raggedright\arraybackslash}p{0.05\textwidth - 2\tabcolsep}>{\raggedright\arraybackslash}p{0.95\textwidth-\tabcolsep}@{}}
\caption{Irreducible diagonal subgroups of $\bar{A}_1^4 A_1^2 = D_7(\#\dseven{17})$. \label{d7diag17}} \\
\hline
ID & Diagonal embedding \\
\hline

\dseven{47} & $(a,a^{[r]},b,c,d,e)$ (if $p=2$ then $r \neq 0$) \\
\dseven{48} & $(a,b,a^{[r]},c,d,e)$ \\
\dseven{49} & $(a^{[r]},b,c,d,a^{[s]},d)$ $(rs=0)$ \\
\dseven{50} & $(a,b,c,d,e,e^{[r]})$ $(r \neq 0)$ \\
\dseven{51} & $(a^{[r]},a^{[s]},a^{[t]},b,c,d)$ ($rt=0; r \leq s;$ if $p=2$ then $r < s$) \\
\dseven{52} & $(a^{[r]},a^{[s]},b,c,a^{[t]},d)$ ($rt=0; r \leq s;$ if $r = s$ then $r < t$ and $p \neq 2$) \\
\hline
\end{longtable}
\vspace{-0.1cm}
\begin{longtable}{llll}
\caption{Irreducible diagonal subgroups of $\bar{A}_1^2 A_1^2 = D_7(\#\dseven{85})$ $(p \neq 2)$. \label{d7diag85}} \\
\hline
ID & Diagonal embedding & ID & Diagonal embedding \\
\hline

\dseven{88} & $(a,a^{[r]},b,c)$ $(r \neq 0)$ & \dseven{92} & $(a^{[r]},b,a^{[s]},a^{[t]})$ $(rs=0; s < t)$ \\

\dseven{89} & $(a^{[r]},b,a^{[s]},c)$ $(rs=0)$ & \dseven{93} & $(a,a^{[r]},b,b^{[s]})$ $(rs \neq 0)$ \\

\dseven{90} & $(a,b,c,c^{[r]})$ $(r \neq 0)$ & \dseven{94} & $(a^{[r]},b^{[t]},a^{[s]},b^{[u]})$ $(rs=tu=0; r \leq t;$ if $r=t$ then $s \leq u$) \\

\dseven{91} & $(a^{[r]},a^{[s]},a^{[t]},b)$ $(rt=0; r < s)$ & \dseven{95} & $(a^{[r]},a^{[s]},a^{[t]},a^{[u]})$ $(rt=0; r < s; t < u)$ \\

\hline
\end{longtable}
\vspace{-0.1cm}
\begin{longtable}{llll}
\caption{Irreducible diagonal subgroups of $\bar{A}_1^2 A_1 A_1 = D_7(\#\dseven{105})$ $(p \geq 7)$. \label{d7diag105}} \\
\hline
ID & Diagonal embedding & ID & Diagonal embedding \\
\hline

\dseven{112} & $(a,a^{[r]},b,c)$ & \dseven{117} & $(a^{[r]},a^{[s]},b,a^{[t]})$ $(rst=0)$  \\
\dseven{113} & $(a^{[r]},b,a^{[s]},c)$ $(rs=0)$ & \dseven{118} & $(a^{[r]},b,a^{[s]},a^{[t]})$ $(rst=0)$  \\
\dseven{114} & $(a^{[r]},b,c,a^{[s]})$ $(rs=0)$ & \dseven{119} & $(a,a^{[r]},b^{[s]},b^{[t]})$ $(st=0)$  \\
\dseven{115} & $(a,b,c^{[r]},c^{[s]})$ $(rs=0)$ & \dseven{120} & $(a^{[r]},b^{[t]},a^{[s]},b^{[u]})$ $(rs=tu=0)$ \\  
\dseven{116} & $(a^{[r]},a^{[s]},a^{[t]},b)$ $(rt=0; r \leq s;$ if $r=s$ then $r < t$) & \dseven{121} & $(a^{[r]},a^{[s]},a^{[t]},a^{[u]})$ $(rtu=0; r \leq s$; \\
& & & if $r=s$ then $r < t$) \\
\hline
\end{longtable}
\vspace{-0.1cm}
\begin{longtable}{ll}
\caption{Irreducible diagonal subgroups of $\bar{A}_1^2 A_1^2 B_2 = D_7(\#\dseven{100})$ $(p=2)$. \label{d7diag100}} \\
\hline
ID & Diagonal embedding \\
\hline

\dseven{123} & $(a,a^{[r]},b,c,10)$ $(r \neq 0)$ \\
\dseven{124} & $(a^{[r]},b,a^{[s]},c,10)$ $(rs=0)$   \\
\dseven{125} & $(a,b,c,c^{[r]},10)$ $(r \neq 0)$   \\
\dseven{126} & $(a^{[r]},a^{[s]},a^{[t]},b,10)$ $(rt=0; r < s)$ \\
\dseven{127} & $(a^{[r]},b,a^{[s]},a^{[t]},10)$ $(rs=0; s < t)$ \\
\dseven{128} & $(a,a^{[r]},b,b^{[s]},10)$ $(rs \neq 0)$ \\
\dseven{129} & $(a^{[r]},b^{[t]},a^{[s]},b^{[u]},10)$ $(rs=tu=0; r \leq t$; if $r=t$ then $s \leq u$)  \\
\dseven{130} & $(a^{[r]},a^{[s]},a^{[t]},a^{[u]},10)$ $(rt=0; r < s; t < u)$ \\
\hline
\end{longtable}
\vspace{-0.1cm}
\begin{longtable}{ll}
\caption{Irreducible diagonal subgroups of $\bar{A}_1^2 A_1^4 = D_7(\#\dseven{122})$ $(p=2)$. \label{d7diag122}} \\
\hline
ID & Diagonal embedding \\
\hline

\dseven{131} & $(a,a^{[r]},b,c,d,e)$ $(r \neq 0)$ \\
\dseven{132} & $(a^{[r]},b,a^{[s]},c,d,e)$ $(rs =0)$ \\
\dseven{133} & $(a,b,c,c^{[r]},d,e)$ $(r \neq 0)$ \\
\dseven{134} & $(a^{[r]},a^{[s]},a^{[t]},b,c,d)$ $(rt=0; r<s)$ \\
\dseven{135} & $(a^{[r]},b,a^{[s]},a^{[t]},c,d)$ $(rs=0; s<t)$ \\
\dseven{136} & $(a,b,c,c^{[r]},c^{[s]},d)$ $(0 < r < s)$ \\
\hline
\end{longtable}

\begin{longtable}{llll}
\caption{Irreducible diagonal subgroups of $\bar{A}_1^2 B_2^2 = D_7(\#\dseven{13})$. \label{d7diag13}} \\
\hline
ID & Diagonal embedding & ID & Diagonal embedding \\
\hline

\dseven{157} & $(1_a,1_a^{[r]},10_b,10_c)$ (if $p=2$ then $r \neq 0)$ & \dseven{160} & $(1_a,1_a^{[r]},10_b,10_b^{[s]})$ ($s \neq 0$; if $p=2$ then $r \neq 0)$  \\

\dseven{158} & $(1_a,1_b,10_c,10_c^{[r]})$ ($r \neq 0)$ & \dseven{161} & $(1_a,1_a^{[r]},10_b,02_b^{[s]})$ ($r \neq 0)$ \\

\dseven{159} & $(1_a,1_b,10_c,02_c^{[r]})$ & & \\

\hline
\end{longtable}

\begin{longtable}{llll}
\caption{Irreducible diagonal subgroups of $\bar{A}_1^2 A_1^2 = D_7(\#\dseven{162})$ $(p \geq 5)$. \label{d7diag162}} \\
\hline
ID & Diagonal embedding & ID & Diagonal embedding \\
\hline

\dseven{166} & $(a,a^{[r]},b,c)$ & \dseven{170} & $(a^{[r]},b,a^{[s]},a^{[t]})$ $(rs=0; s < t)$ \\

\dseven{167} & $(a^{[r]},b,a^{[s]},c)$ $(rs=0)$ & \dseven{171} & $(a,a^{[r]},b,b^{[s]})$ $(s \neq 0)$ \\

\dseven{168} & $(a,b,c,c^{[r]})$ $(r \neq 0)$ & \dseven{172} & $(a^{[r]},b^{[t]},a^{[s]},b^{[u]})$ $(rs=tu=0; r \leq t;$ if $r=t$ then $s \leq u$) \\

\dseven{169} & $(a^{[r]},a^{[s]},a^{[t]},b)$ $(rt=0; r \leq s)$ & \dseven{173} & $(a^{[r]},a^{[s]},a^{[t]},a^{[u]})$ $(rt=0; r < s; t < u)$ \\
\endfirsthead
\hline
\end{longtable}

\begin{longtable}{llll}
\caption{Irreducible diagonal subgroups of $A_1 A_1 A_1 = D_7(\#\dseven{189})$ $(p \geq 5)$. \label{d7diag189}} \\
\hline
ID & Diagonal embedding & ID & Diagonal embedding \\
\hline

\dseven{191} & $(a^{[r]},a^{[s]},b)$ $(rs=0)$ & $\mathit{\dseven{169}^{\{\underline{0}\}}}$ & $(a,b,b)$ \\
\dseven{192} & $(a^{[r]},b,a^{[s]})$ $(rs=0)$ & \dseven{194} & $(a^{[r]},a^{[s]},a^{[t]})$ $(rst=0; s \neq t)$ \\
\dseven{193} & $(a,b^{[r]},b^{[s]})$ $(rs=0; r \neq s)$ \\
\hline
\end{longtable}

\begin{longtable}{ll}
\caption{Irreducible diagonal subgroups of $A_1^2 A_1^3 = D_7(\#\dseven{200})$ $(p = 2)$. \label{d7diag200}} \\
\hline
ID & Diagonal embedding \\
\hline

\dseven{203} & $(a,a^{[r]},b,c,d)$ $(r \neq 0)$ \\
\dseven{204} & $(a^{[r]},b,a^{[s]},c,d)$ $(rs=0)$ \\
\dseven{205} & $(a,b,c,c^{[r]},d)$ $(r \neq 0)$ \\
\endfirsthead
\dseven{206} & $(a^{[r]},a^{[s]},a^{[t]},b,c)$ $(rt=0; r < s)$ \\
\dseven{207} & $(a^{[r]},b,a^{[s]},a^{[t]},c)$ $(rs=0; s < t)$ \\
\dseven{208} & $(a,b,c,c^{[r]},c^{[s]})$ $(0 < r < s)$ \\
\dseven{209} & $(a,a^{[r]},b,b^{[s]},c)$ $(rs \neq 0)$ \\
\dseven{210} & $(a^{[r]},b^{[t]},a^{[s]},b^{[u]},c)$ $(rs=tu=0; r \leq t;$ if $r=t$ then $s \leq u$) \\
\dseven{211} & $(a^{[r]},b,a^{[s]},c,c^{[t]})$ $(rs=0 \neq t)$ \\
\dseven{212} & $(a^{[r]},a^{[s]},a^{[t]},a^{[u]},b)$ $(rt=0; r < s; t < u)$ \\
\dseven{213} & $(a^{[r]},b,a^{[s]},a^{[t]},a^{[u]})$ $(rs=0; s < t < u)$ \\
\dseven{214} & $(a^{[r]},a^{[s]},a^{[t]},b,b^{[u]})$ $(rt=0 \neq u; r < s)$ \\
\dseven{215} & $(a^{[r]},b^{[u]},a^{[s]},a^{[t]},b^{[v]})$ $(rs=uv=0; s < t)$ \\
\dseven{216} & $(a,a^{[r]},b,b^{[s]},b^{[t]})$ $(r \neq 0; 0 < s < t)$ \\
\dseven{217} & $(a^{[r]},a^{[s]},a^{[t]},a^{[u]},a^{[v]})$ $(rt=0; r < s; t < u < v)$ \\
\hline
\end{longtable}

\begin{longtable}{llll}
\caption{Irreducible diagonal subgroups of $A_1^3 A_1 = D_7(\#\dseven{224})$ $(p \geq 5)$. \label{d7diag224}} \\
\hline
ID & Diagonal embedding & ID & Diagonal embedding \\
\hline
\dseven{227} & $(a,a^{[r]},b,c)$ $(r \neq 0)$ & \dseven{230} & $(a^{[r]},a^{[s]},b,a^{[t]})$ $(rt=0; r < s)$ \\
\dseven{228} & $(a^{[r]},b,c,a^{[s]})$ $(rs=0)$ & \dseven{231} & $(a,a^{[r]},b^{[s]},b^{[t]})$ $(r \neq 0 = st)$ \\
\dseven{229} & $(a,a^{[r]},a^{[s]},b)$ $(0 < r < s)$ & \dseven{232} & $(a^{[r]},a^{[s]},a^{[t]},a^{[u]})$ $(ru=0; r < s < t)$ \\
\hline \endfirsthead
\end{longtable}

\begin{longtable}{llll}
\caption{Irreducible diagonal subgroups of $A_1^2 A_1 A_1 = D_7(\#\dseven{233})$ $(p \geq 5)$. \label{d7diag233}} \\
\hline
ID & Diagonal embedding \\
\hline
\dseven{237} & $(a,a^{[r]},b,c)$ $(r \neq 0)$ & \dseven{242} & $(a^{[r]},a^{[s]},b,a^{[t]})$ $(rt=0; r < s)$ \\
\dseven{238} & $(a^{[r]},b,a^{[s]},c)$ $(rs=0)$ & \dseven{243} & $(a^{[r]},b,a^{[s]},a^{[t]})$ $(rst=0; s \neq t)$ \\
\dseven{239} & $(a^{[r]},b,c,a^{[s]})$ $(rs=0)$ & \dseven{244} & $(a,a^{[r]},b^{[s]},b^{[t]})$ $(r \neq 0 = st; s \neq t)$ \\
\endfirsthead
\dseven{240} & $(a,b,c^{[r]},c^{[s]})$ $(rs = 0; r \neq s)$ & \dseven{245} & $(a^{[r]},b^{[t]},a^{[s]},b^{[u]})$ $(rs=tu=0)$ \\
$\mathit{\dseven{228}^{\{\underline{0}\}}}$ & $(a,b,c,c)$  & \dseven{246} & $(a^{[r]},a^{[s]},a^{[t]},a^{[u]})$ $(rtu=0; r < s; t \neq u)$ \\
\dseven{241} & $(a^{[r]},a^{[s]},a^{[t]},b)$ $(rt=0; r < s)$ \\
\hline
\end{longtable}

\begin{longtable}{llll}
\caption{Irreducible diagonal subgroups of $A_1^6 = D_7(\#\dseven{283})$ $(p =2)$. \label{d7diag283}} \\
\hline
ID & Diagonal embedding \\
\hline

\dseven{288} & $(a,a^{[r]},b,c,d,e)$ $(r \neq 0)$ & \dseven{294} & $(a,a^{[r]},a^{[s]},a^{[t]},a^{[u]},b)$ $(0 < r < s < t < u)$ \\
\dseven{289} & $(a,a^{[r]},a^{[s]},b,c,d)$ $(0 < r < s)$ & \dseven{295} & $(a,a^{[r]},a^{[s]},a^{[t]},b,b^{[u]})$ $(0 < r < s < t; u \neq 0)$ \\
\dseven{290} & $(a,a^{[r]},b,b^{[s]},c,d)$ $(rs \neq 0; r \leq s)$ & \dseven{296} & $(a,a^{[r]},a^{[s]},b,b^{[t]},b^{[u]})$ $(0 < r < s; 0 < t < u; r \leq t;$ \\
\dseven{291} & $(a,a^{[r]},a^{[s]},a^{[t]},b,c)$ $(0 < r < s < t)$ & & if $r =t$ then $s \leq u$) \\ 
\dseven{292} & $(a,a^{[r]},a^{[s]},b,b^{[t]},c)$ $(0 < r < s; t \neq 0)$ & \dseven{297} & $(a,a^{[r]},a^{[s]},a^{[t]},a^{[u]},a^{[v]})$ $(0 < r < s < t < u < v)$ \\ 

\dseven{293} & $(a,a^{[r]},b,b^{[s]},c,c^{[t]})$ $(rst \neq 0; r \leq s \leq t)$ & & \\
\hline
\end{longtable}

\begin{longtable}{llll}
\caption{Irreducible diagonal subgroups of $A_1^2 B_2^2 = D_7(\#\dseven{279})$ $(p =2)$. \label{d7diag279}} \\
\hline
ID & Diagonal embedding \\
\hline
\dseven{298} & $(1_a,1_b^{[r]},10_c,10_d)$ ($r \neq 0)$ & \dseven{301} & $(1_a,1_a^{[r]},10_b,10_b^{[s]})$ ($rs \neq 0$) \\
\dseven{299} & $(1_a,1_b,10_c,10_c^{[r]})$ ($r \neq 0)$ & \dseven{302} & $(1_a,1_a^{[r]},10_b,02_b^{[s]})$ ($r \neq 0$) \\
\endfirsthead
\dseven{300} & $(1_a,1_b,10_c,02_c^{[r]})$ \\
\hline
\end{longtable}
}

\section{Applying the main theorem} \label{sec:applying}

We now describe how Theorem \ref{thm:main} allows one to classify all connected $G$-cr subgroups, and illustrate with the case $G = F_{4}(K)$, in arbitrary characteristic. In this section we present a series of additional lemmas which allow one to derive properties of the $G$-cr subgroups arising. The subgroups themselves, and various properties thereof, are given in Table \ref{tab:f4}.

\subsection{$G$-conjugacy vs.\ $L$-conjugacy} \label{sec:conjugacy}

Once we know the $L'$-irreducible subgroups for each possible Levi subgroup $L$ of our reductive algebraic group $G$, it remains to consider how the various classes of subgroups fuse in $G$. This is achieved in the following result.

\begin{proposition} \label{prop:glconj}
Let $X$ and $Y$ be $G$-cr subgroups of $G$, and let $L$ and $L^{\ast}$ be Levi subgroups of $G$ which are minimal among those containing $X$ and $Y$ respectively. Suppose further that $X$ and $Y$ are $G$-conjugate. Then $L$ and $L^{\ast}$ are also $G$-conjugate; moreover if $L^{\ast} = L$ then $X$ and $Y$ are in fact $N_{G}(L)$-conjugate.
\end{proposition}

\proof Take $g \in G$ such that $Y = X^{g}$. Then $Z(L)^{\circ}$ and $(Z(L^{\ast})^{\circ})^{g^{-1}}$ are both maximal tori of $C_{G}(X)^{\circ}$, so there exists $c \in C_{G}(X)^{\circ}$ such that $(Z(L^{\ast})^{\circ})^{g^{-1}} = (Z(L)^{\circ})^{c}$. Thus $L^{\ast} = C_{G}(Z(L^{\ast})^{\circ}) = C_{G}(Z(L)^{\circ})^{cg} = L^{cg}$, proving the first statement. Now assume $L^{\ast} = L$, and write $n = cg$. Then $L^{n} = L$ so $n \in N_{G}(L)$, and $X^{n} = X^{cg} = X^{g} = Y$, proving the second statement. \qed

In light of the above result, it is desirable to understand how the classes of subgroups of simple Levi factors $L_{0}$ fuse under graph automorphisms of $L_{0}$. We consider only the case that $L_{0}$ has type $A_n$ or $D_{n}$ $(n \ge 1)$, since the exceptional graph morphism of $B_{2}$ $(p = 2)$ is not a variety automorphism. For $L_{0}$ of type $A_{n}$, we refer to the well-known fact that an irreducible module for a connected reductive algebraic group is either self-dual, or equivalent to its dual under an outer automorphism. It follows that every irreducible reductive subgroup of $L_{0}$ is normalised by an outer automorphism of $L_{0}$, hence by Lemma \ref{lem:go_so_classes} distinct $L_{0}$-conjugacy classes of such subgroups are not fused by an outer automorphism.

It remains to consider $L_{0}$ of type $D_n$. Now, $D_7$ only occurs as a Levi factor in $E_8$, and its normaliser then contains an element inducing a graph automorphism, and the irreducible subgroups in Table \ref{D6tab} are already given up to $\Aut(D_7)$-conjugacy. This leaves $L_{0}$ of type $D_4$, $D_5$ and $D_6$. In these cases, the tables present the irreducible subgroups up to $L_{0}$-conjugacy, and we describe how these subgroup classes fuse under graph automorphisms. Firstly, for types $D_5$ and $D_6$ the notation `$l \sim m$' indicates that $D_n(\#l)$ is $\Aut(D_n)$-conjugate to $D_n(\#m)$. For type $D_6$, the notation `$l:$ condition' for the subgroup $D_n(\#l)$ indicates that $D_{n}(\#l)$ is a collection of diagonal subgroups and the `condition' is an extra constraint that is placed on the corresponding field twists to obtain the representatives of $\Aut(D_n)$-classes of subgroups. For type $D_4$ we are interested in the fusion of classes under both an involutory graph automorphism and under the full outer automorphism group $S_3$. We present this using brackets, for instance, $\dfour{3} \sim \dfour{4} (\sim \dfour{2})$ indicates that the two $B_3$ classes with ID numbers $\dfour{3}$ and $\dfour{4}$ are fused by an involutory graph automorphism, and that the three $B_3$ classes with ID numbers $\dfour{2},\dfour{3}$ and $\dfour{4}$ are fused under the full outer automorphism group. 

$D_4$: $\dfour{3} \sim \dfour{4}$ $(\sim \dfour{2})$; $\dfour{6} \sim \dfour{7}$ $(\sim \dfour{5})$; $\dfour{10} \sim \dfour{11}$ $(\sim \dfour{9})$; $\dfour{12}$: $r \leq s$ $(s \leq t)$; $\dfour{14} \sim \dfour{15}$ $(\sim \dfour{13})$; $\dfour{16}$: $s \leq t$ ($r \leq s$); $\dfour{24} \sim \dfour{29}$ $(\sim \dfour{19})$; $\dfour{25} \sim \dfour{30}$ $(\sim \dfour{20})$; $\dfour{26} \sim \dfour{31}$ $(\sim \dfour{21})$; $\dfour{27} \sim \dfour{32}$ $(\sim \dfour{22})$; $\dfour{28} \sim \dfour{33}$ $(\sim \dfour{23})$.  

$D_5$: $\dfive{5} \sim \dfive{6}$

$D_6$: $\dsix{6} \sim \dsix{7}$; $\dsix{18} \sim \dsix{19}$; $\dsix{21}$: $r \leq s$; $\dsix{23}$: $r \leq t$, if $r=t$ then $s \leq u$; $\dsix{25}$: $r \leq s$; $\dsix{26}$: $r \leq s$; $\dsix{27}$: $r \leq s$; $\dsix{28} \sim \dsix{29}$; $\dsix{30}$: $t \leq u$; $\dsix{31}$: $r \leq s$; $\dsix{32}$: $r \leq t \leq v$, if $r=t$ then $s \leq u$, if $t=v$ then $u \leq w$; $\dsix{33}$: $t \leq u$; $\dsix{34}$: $s \leq t$; $\dsix{35}$: $r \leq s$; $\dsix{36}$: $v \leq w$; $\dsix{37}$: if $t=0$ then $r \leq u$, if $t=0$ and $r=u$ then $s \leq v$; $\dsix{38}$: in Table \ref{condD639} remove line 6 and add the following to the further requirements column: $v < w$ to lines 1 and 3, $u<w$ to line 8, $r<s<u$ to line 11; $\dsix{42}$: $r \leq s$; $\dsix{75}$: $r \leq s$; $\dsix{76} \sim \dsix{77}$; $\dsix{78}$: $r \leq s$; $\dsix{79}$: $r \leq s$; $\dsix{80} \sim \dsix{82}$; $\dsix{81} \sim \dsix{83}$; $\dsix{85}$: $r \leq s$; $\dsix{86}$: $r \leq s$; $\dsix{87} \sim \dsix{88}$; $\dsix{89}$: $r \leq s$; $\dsix{90} \sim \dsix{91}$; $\dsix{92}$: $r \leq s$; $\dsix{93}$: $r < s$; $\dsix{94} \sim \dsix{95}$; $\dsix{97}$: $r < s$; $\dsix{98} \sim \dsix{99}$; $\dsix{101}$: $r < s$; $\dsix{102}$: $r \leq t$ and if $r=t$ then $s \leq u$; $\dsix{103} \sim \dsix{104}$; $\dsix{105}$: $r < s$; $\dsix{106} \sim \dsix{107}$; $\dsix{108}$: $r < s$; $\dsix{109} \sim \dsix{110}$; $\dsix{111}$: $r < s$; $\dsix{112}$: $r < s$; $\dsix{116} \sim \dsix{117}$; $\dsix{155} \sim \dsix{163}$; $\dsix{156} \sim \dsix{164}$; $\dsix{157} \sim \dsix{165}$; $\dsix{158} \sim \dsix{160}$; $\dsix{159}$: $r < s$; $\dsix{162} \sim \dsix{166}$.

\subsection{Normaliser Structure}

Let $X$ be a $G$-cr subgroup of $G$. The following result gives us a method for calculating the structure of $N_{G}(X)$. In particular, one need only inspect the parabolic subgroups of $G$ which are minimal subject to containing $X$, and the normaliser of a Levi subgroup which is minimal subject to containing $X$. By Proposition \ref{prop:glconj} there is a unique such Levi subgroup up to conjugacy, and there are finitely many such parabolic subgroups to consider, corresponding to the different standard parabolic subgroups having conjugate Levi factors.\\

\begin{proposition} \label{prop:normaliser}
Let $X$ be a $G$-cr subgroup of $G$.
\begin{enumerate}[label=\normalfont(\roman*)]
	\item If $L$ is minimal among Levi subgroups of $G$ containing $X$, then
\begin{align*}
N_{G}(X) &= C_{G}(X)^{\circ}(N_{G}(L) \cap N_{G}(X)).
\end{align*}
	\item There exists a parabolic subgroup $P$, with Levi decomposition $P = QL$, such that $P$ is minimal among parabolic subgroups containing $X$ and $C_{G}(X)^{\circ} = \left< C_{Q}(X), C_{Q^{\rm op}}(X),Z(L)^{\circ} \right>$, where $Q^{\rm op}$ is the unipotent radical of the parabolic subgroup opposite to $P$.
\end{enumerate}
\end{proposition}

\proof (i) Taking $Y = X$ in the proof of Proposition~\ref{prop:glconj}, we see that if $g \in N_{G}(X)$ then $g = c^{-1}n$ for some $c \in C_{G}(X)^{\circ}$ and some $n \in N_{G}(L) \cap N_{G}(X)$.

(ii) Let $U$ be a maximal unipotent subgroup of $C_{G}(X)^{\circ}$. Then there is a parabolic subgroup $P$ of $G$ which contains $UX$, such that $U \le Q \deq R_{u}(P)$. Since $X$ is $G$-cr, there exists a Levi subgroup $L$ of $P$ containing $X$, and moreover we can assume that $X$ is $L$-irreducible, otherwise there is a smaller parabolic subgroup of $G$ containing both $UX$ and the unipotent radical of a proper parabolic subgroup of $L$ normalised by $X$. Then $U \le C_{Q}(X)$, and thus $U = C_{Q}(X)$ by maximality of $U$. Since $X$ is $L$-irreducible, $Z(L)^{\circ}$ is a maximal torus of $C_{G}(X)^{\circ}$. Therefore $UZ(L)^{\circ}$, which contains a maximal torus of $C_{G}(X)^{\circ}$ and a maximal unipotent subgroup of $C_{G}(X)^{\circ}$, is a Borel subgroup of $C_{G}(X)^{\circ}$. For the same reason, $C_{Q^{\rm op}}(X)Z(L)^{\circ}$ is also a Borel subgroup of $C_{G}(X)^{\circ}$, and is clearly opposite to $C_{Q}(X)Z(L)^{\circ}$. Thus $C_{G}(X)^{\circ}$ is generated by $C_{Q}(X)Z(L)^{\circ}$ and $C_{Q^{\rm op}}(X)Z(L)^{\circ}$, which gives the desired result. \qed

\subsection{Action on $G$-modules, separability and reductive pairs} \label{sec:action_sep_rp}

Given a reductive subgroup $X$ of $G$, it is of interest to know how $X$ acts on various $G$-modules, particularly the Lie algebra $L(G)$ and the non-trivial $G$-module(s) of minimal dimension. Such information allows one, for instance, to study the conjugacy classes of $G$ meeting $X$. Additionally, recall that $X$ is called \emph{separable} in $G$ if $\dim C_{G}(X) = \dim C_{L(G)}(X)$; equivalently, the scheme-theoretic centraliser of $X$ in $G$ is smooth. Once the normaliser structure and action of $X$ on $L(G)$ are known, one can see directly whether this equality holds. Recall also that $(G,X)$ is a \emph{reductive pair} if $L(X)$ is an $X$-module direct summand of $L(G)$. Both of these properties are closely related to complete reducibility, cf.\ \cite[Theorem 3.35]{MR2178661}.

The next lemma is useful in determining the action of $X$ on $L(G)$. Specifically, one first considers the action of $L$, where $L$ is minimal among Levi subgroups of $G$ containing $X$, and then considers the action of $X$ on each $L$-module arising. We briefly recall some material from \cite{MR1047327}. For a subset $I$ of simple roots of $G$, each root $\beta$ of $G$ can be written uniquely as $\beta_{I} + \beta_{I'}$ where $\beta_{I}$ and $\beta_{I'}$ are respectively a linear combination of the simple roots in $I$ and the simple roots not in $I$. The \emph{shape} of $\beta$ (with respect to $I$) is defined to be $\beta_{I'}$. For a fixed shape $S$, we denote by $V_{S}$ the sum of all root subspaces of $L(G)$ corresponding to roots of shape $S$. Since each root subgroup of $G$ is isomorphic as a $T$-module to its Lie algebra, this definition of $V_{S}$ is compatible with that given in \cite[p.\ 553]{MR1047327}.

\begin{lemma} \label{lem:action}
Let $G$ be a reductive algebraic group and let $L$ be a Levi subgroup of $G$.
\begin{enumerate}[label=\normalfont(\roman*)]
\item If $V$ is a tilting module for $G$ then it is tilting for $L$ and for $L'$;
\item Suppose $L = L_{I}$ is a standard Levi subgroup, corresponding to some subset $I$ of the simple roots of $G$. Then as $L$-modules, we have a direct sum decomposition
\[ L(G) \downarrow L = L(L) \oplus \bigoplus V_{S}, \]
the right-hand sum taken over non-zero shapes $S$.
\end{enumerate}
\end{lemma}

\proof (i) is \cite[Proposition 1.2(ii), Lemma 1.4(i)]{DonkinTilting}. (ii) $L$ is generated by some maximal torus $T$ of $G$ together with those root subgroups corresponding to roots of shape zero (i.e.\ the roots in the span of $I$). Similarly $L(L)$ is generated by $L(T)$ and the root subspaces corresponding to roots of shape zero. Hence as vector spaces, the given direct sum holds. Since $T$ preserves each root subspace, and since a root subgroup of $L$ acting on a root vector of non-zero shape $S$ can only produce sums of other root vectors of shape $S$, it follows that each $V_{S}$ in the given direct sum is an $L$-submodule of $L(G)$ (cf.\ \cite[p.\ 64]{MR0407163} for the action of root elements on root spaces). \qed

Importantly, the modules $V_{S}$ arising in \cite{MR1047327} have a very limited range of possible high weights, which can be explicitly determined using \cite[Theorem 2 and Remark 1]{MR1047327}. For $G$ of exceptional type, the weights occurring are given in \cite[Lemma 3.1]{MR1329942}. In most cases, these modules lie in the tensor algebra of the natural modules for the simple classical factors of $L$.

\subsection{Subgroups of $F_{4}$} \label{sec:f4table}

We now classify the $G$-completely reducible semisimple subgroups of $G = F_{4}(K)$. For each such subgroup $X$, we also give the $X$-module structure of $L(G)$, as well as the connected centraliser $C^{\circ} \deq C_{G}(X)^{\circ}$. Using these, we also determine whether $X$ is a separable subgroup of $G$ (`Sep'), and whether $(G,X)$ is a reductive pair (`RP'). More precisely, in the two rightmost columns of Table \ref{tab:f4} below either we write `Yes' or `No' to indicate whether the subgroup is separable or forms a reductive pair with $G$, or we write conditions on the characteristic and field twists to indicate that the subgroup is separable (or forms a reductive pair with $G$) precisely when these conditions hold.

Our process for classifying subgroups is as follows. For each Levi subgroup $L$, and each simple factor $L_{0}$ of $L$, the image in $L_{0}$ of an $L$-irreducible subgroup $X$ of $L$ is $L_{0}$-conjugate to one of the subgroups in Section \ref{sec:tables}. Thus each simple factor of $X$ either is contained in a simple factor of $L$, or is a diagonal subgroup of a product of two or more such factors.

Let us illustrate with the case $L_{134}$, of type $A_{1} A_{2}$. Then $X$ projects surjectively to the first factor, and its image in the second factor is either $A_2$ itself or an irreducible subgroup of type $A_1$ (given by $A_{2}(\#1)$ in Table \ref{A2tab}). Thus one of the following holds:
\begin{itemize}
	\item $X = L'$. We denote this with the ID $L_{134}(\#0)$;
	\item $p \neq 2$ and $X = A_{1}A_{1}$, where the second factor is the subgroup $A_{2}(\#1)$ in Table \ref{A2tab}. We denote this by $L_{134}(\#0;\#1)$;
	\item $p \neq 2$, $X = A_{1}$, and $X$ is embedded diagonally in the subgroup $L_{134}(\#0;\#1)$ above. Such an embedding is determined by a pair of Frobenius twists, and we therefore denote this by $L_{134}(\#0^r,\#1^s)$, where $r$ and $s$ are the powers of $p$ giving the twists. In such a situation we implicitly take $rs = 0$.
\end{itemize}

Next, the stated composition factors of $X = L_{134}(\#0) = L_{134}'$ on $L(G)$ follow from Lemma \ref{lem:action} and the calculation of the high weights of $L$ on the modules $V_{S}$. When $p \neq 2$, $L(G)$ is tilting, and we obtain the direct sum decomposition stated. When $p = 2$, we consider instead the tilting modules $T_{G}(\lambda_4) = \lambda_4$ (irreducible of dimension 26) and $T_{G}(\lambda_1) = \lambda_4 \mid \lambda_1 \mid \lambda_4$. Each of these restricts to $X$ as a tilting module. Thus, using for instance Doty's software \cite{Dot1} for computing the structure of Weyl modules, we deduce that
\begin{align*}
\begin{split}
T(\lambda_1) \downarrow X &= (1,T(20)) + (1,T(02)) + (0,11)^{2} \\
& \qquad + (0,T(20)) + (0,T(02)) + (T(2),0) + (1,0)^{2},\\
\lambda_4 \downarrow X &= (1,01) + (1,10) + (0,01) + (0,10) + (0,11),
\end{split}
\end{align*}
where we have $T(2) = 0 \mid 2 \mid 0$ for $A_1$, and $T(20) = 01 \mid 20 \mid 01$ for $A_{2}$. Therefore $L(G)$, a $52$-dimensional submodule of $T(\lambda_1)$ with $\lambda_4$ as a submodule, has structure
\begin{align*}
\begin{split}
L(G) \downarrow X &= (1,W(20)) + (0,W(20)) + (0,11)\\
& \qquad + (1,W(02)) + (0,W(02)) + (T(2),0) + (1,0)^2
\end{split}
\end{align*}
as in the table. For the other possibilities for $X$, we need only compute the restriction of the modules for $L_{134}'$, which are tensor products of modules for $X$, to deduce the given module structure.

Finally, each possible subgroup $X$ above clearly centralises the $1$-dimensional torus $Z(L_{134})$. Since the centraliser of a $G$-cr subgroup is $G$-cr \cite[Corollary 3.17]{MR2178661} and thus reductive, we deduce that $C_{G}(X)^{\circ}$ is either a $1$-dimensional torus (`$T_1$') or a simple subgroup of type $A_{1}$. Using Proposition \ref{prop:normaliser}, we need to consider the action of $X$ on the unipotent radical $Q$ of the (unique) standard parabolic subgroup having $L_{134}$ as a Levi factor. We find that the proper subgroups of $L_{134}'$ (and only these) centralise a vector of $Q$. Specifically, the module $V_{S}$, where $S$ consists of roots having $\alpha_2$-coefficient equal to $2$, is an irreducible high weight module $(0,20)$ for $L_{134}'$, and thus restricts to the diagonal subgroup $L_{134}(\#0;\#1)$ as the symmetric square $S^{2}(2) = 4 + 0$. The $1$-dimensional trivial submodule here gives rise to a $1$-dimensional $X$-invariant subgroup, hence $\dim C_{G}(X)^{\circ} \ge 2$, and so $C_{G}(X)^{\circ}$ is simple of type $A_{1}$. This gives all the information for these subgroups in Table \ref{tab:f4}.\\

\begin{remark} \label{rem:f4}
The following semisimple reducible $G$-cr subgroups were omitted in \cite[Table 4.1]{MR3075783}:
$L_{124}(\#1;\#0)$, $L_{124}(\#1^{[r]};\#0^{[s]})$, $L_{134}(\#0;\#1)$, $L_{134}(\#0^{[r]};\#1^{[s]})$, $L_{123}(\#5)$, $L_{123}(\#6)$, $L_{123}(\#7)$ when $p=2$, $L_{123}(\#11)$, $L_{123}(\#12)$, $L_{234}(\#9)$, $L_{234}(\#10)$, $L_{234}(\#11)$, $L_{234}(\#13)$, $L_{234}(\#14)$, $L_{234}(\#15)$ when $p \geq 3$, $L_{23}(\#3)$ when $p=2$, $L_{23}(\#5)$, and $L_{34}(\#1)$ when $p \geq3$. 
\end{remark}

{ \small
\begin{longtable}{
>{\raggedright\arraybackslash}p{0.07\textwidth - 2\tabcolsep} 
>{\raggedright\arraybackslash}p{0.16\textwidth-2\tabcolsep} 
>{\raggedright\arraybackslash}p{0.16\textwidth-2\tabcolsep}
>{\raggedright\arraybackslash}p{0.395\textwidth-2\tabcolsep} 
>{\raggedright\arraybackslash}p{0.08\textwidth-2\tabcolsep} 
>{\raggedright\arraybackslash}p{0.06\textwidth-2\tabcolsep} 
>{\raggedright\arraybackslash}p{0.08\textwidth-\tabcolsep}@{} 
}
\caption{Reducible, $G$-cr subgroups of $G = F_4(K)$ \label{tab:f4}} \\
\hline
 
$X$ & ID & Conditions & $L(G) \downarrow X$ & $C^{\circ}$ & Sep & RP \\ \hline
$B_3$				& $L_{123}(\#0)$ & $p \neq 2$ & $010 + 100^{2} + 001^{2} + 0$ & $T_{1}$ & Yes & Yes \\
					& & $p = 2$ & $((010 + 0) | 100 | 0) + W(100)^{2} + 001^{2}$ & $\tilde{A}_{1}$ & Yes & No \\ \hline
\endfirsthead
$C_3$				& $L_{234}(\#0)$ & $p \neq 2$ & $200 + 001^{2} + 0^{3}$ & $\bar{A}_{1}$ & Yes & Yes \\
					& & $p = 2$ & $(0 | 200 | (010 + 0)) + W(001)^{2} + 0^{2}$ & $\bar{A}_{1}$ & Yes & No \\ \hline
$A_3$				& $L_{123}(\#2)$ & & $T(101;0) + 100^2 + 010^{3} + 001^2 + 0^{2}$ & $\tilde{A}_{1}$ & Yes & $p \neq 2$ \\
					& $L_{234}(\#3)$ & $p = 2$ & $( (020 + 0) | 101 ) + ((200 + 002) | 010)^{2} + 0^{3}$ & $\bar{A}_{1}$ & Yes & No \\ \hline
$A_{1} B_{2}$		& $L_{123}(\#4)$ & $p = 2$ & $(0|(0,02)|((2,10)+0^{2})|((2,0)+(0,10))|0) + (((2,0)+(0,10)) | 0)^{2} + (1,01)^{2}$ & $\tilde{A}_{1}$ & Yes & No \\
					& $L_{234}(\#1)$ & $p \neq 2$ & $(2,0) + (0,02) + (0,01)^{2} + (1,01) + (1,10)^{2} + 0^{3}$ & $\bar{A}_{1}$ & Yes & Yes \\
					& & $p = 2$ & $(0|(0,02)|0|((0,10)+(2,0))|0^{2}) + (0,01)^{2} + (1,01) + (1,W(10))^{2} + 0^{2}$ & $\bar{A}_{1}$ & No & No \\ \hline
$A_{1} A_{2}$		& $L_{124}(\#0)$ & $p \neq 2$ & $(2,0) + (0,T(11;0)) + (0,10) + (0,01) + (1,10) + (1,01) + (2,10) + (2,01) + (1,0)^{2}$ & $T_1$ & Yes & $p \neq 3$ \\
					& & $p = 2$ & $( T(2),0 ) + (0,11) + (0,10) + (0,01) + (1,10) + (1,01) + (W(2),10) + (W(2),01) + (1,0)^{2}$ & $T_1$ & Yes & No \\
					& $L_{134}(\#0)$ & $p \neq 2$ & $(1,20) + (1,02) + (0,20) + (0,02) + (0,T(11;0)) + (2,0) + (1,0)^2$ & $T_1$ & Yes & $p \neq 3$ \\
					& & $p = 2$ & $(1,W(20)) + (0,W(20)) + (0,11) + (1,W(02)) + (0,W(02)) + (T(2),0) + (1,0)^2$ & $T_1$ & Yes & No \\ \hline
$A_{1}^{3}$	& $L_{123}(\#1)$ & $p \neq 2$ & $(1,1,2) + (1,1,0)^{2} + (0,1,1)^{2} + (1,0,1)^{2} + (2,0,0) + (0,2,0) + (0,0,2)^{3} + 0$ & $T_{1}$ & Yes & Yes \\
					& & $p = 2$    & $(1,1,W(2)) + (1,1,0)^{2} + (0,1,1)^{2} + (1,0,1)^{2} + (0^{2} | (2,0,0)+(0,2,0)+(0,0,2)| 0^{2}) + (0,0,W(2))^{2}$ & $\tilde{A}_{1}$ & No & No \\
					& $L_{123}(\#10)$ & $p = 2$ & $(((2,2,0) + (2,0,2) + (0,2,2) + 0^{3}) | ((2,0,0) + (0,2,0) + (0,0,2)) | 0)+ (((2,0,0)+(0,2,0)+(0,0,2)) | 0)^{2} + (1,1,1)^{2}$ & $\tilde{A}_{1}$ & Yes & No \\
					& $L_{234}(\#6)$ & $p = 2$ & $((0^{2} + (0,2,2)) | ((2,0,0) + (0,2,0) + (0,0,2)) | 0^{2}) + (0,1,1)^{2} + (1,1,1) + (((1,2,0) + (1,0,2))|(1,0,0))^{2} + 0^{2}$ & $\bar{A}_{1}$ & No & No \\
					& $L_{234}(\#8)$ & $p \neq 2$ & $(2,0,0) + (0,2,0) + (0,0,2) + (0,1,1) + (1,1,0) + (1,0,1) + (1,0,0)^{2} + (0,1,0)^{2} + (0,0,1)^{2} + (1,1,1)^{2} + 0^{3}$ & $\bar{A}_{1}$ & Yes & Yes \\
					& & $p = 2$ & $(0 | ((2,0,0) + (0,2,0) + (0,0,2)) | 0^{3}) + (1,1,1)^{2} + (1,1,0) + (0,1,1) + (1,0,1) + (1,0,0)^{2} + (0,1,0)^{2} + (0,0,1)^{2} + 0^{2}$ & $\bar{A}_{1}$ & No & No \\ \hline
$G_{2}$				& $L_{123}(\#3)$ & $p \neq 2$ & $W(01) + 10^{5} + 0^{3}$ & $\tilde{A}_1$ & Yes & Yes \\
					& & $p = 2$ & $01 + W(10)^{2} + T(10)^{3}$ & $\tilde{A}_1$ & No & Yes \\
					& $L_{234}(\#4)$ & $p = 2$ & $(0 | 20 | (01 + 0)) + ( 0 | 20 | 0 | 10 )^{2} + 0^{2}$ & $\bar{A}_{1}$ & Yes & No \\ \hline
$B_{2}$				& $L_{23}(\#0)$ & $p \neq 2$ & $02 + 10^{4} + 01^{4} + 0^{6}$ & $\bar{A}_{1}^2$ & Yes & Yes \\
					& & $p = 2$ & $(0 | 10 | 0 | 02 | 0) + W(10)^{4} + 01^{4} + 0^{5}$ & $B_{2}$ & Yes & No \\ \hline
$A_{2}$				& $L_{12}(\#0)$ & & $T(11;0) + 10^6 + 01^6 + 0^7$ & $\tilde{A}_2$ & Yes & $p \neq 3$ \\
					& $L_{34}(\#0)$ & & $T(11;0) + W(20)^3 + W(02)^3 + 0^7$ & $\bar{A}_2$ & Yes & $p \neq 3$ \\
					& $L_{123}(\#8)$ & $p = 3$ & $((30 + 03 + 0) | 11 ) + 11^{5} + 0^{3}$ & $A_{1}$ & Yes & No \\ \hline
$A_{1}^{2}$		& $L_{13}(\#0)$ & $p \neq 2$ & $(2,0) + (0,2)^{3} + (1,0)^{4} + (0,1)^{4} + (1,1)^{2} + (1,2)^{2} + 0^{4}$ & $\bar{A}_1 T_1$ & Yes & Yes \\
					& & $p = 2$ & $(T(2),0) + (1,W(2))^2 + (1,1)^2 + (1,0)^4 + (0,T(2)) + (0,W(2))^2 + (0,1)^4 + 0^2$ & $\bar{A}_1 \tilde{A}_1$ & Yes & No \\
					& $L_{124}(\#1;\#0)$ & $p \neq 2$ & $(0,T(4;0)) + (2,0) + (0,2)^{3} + (1,2)^{2} + (2,2)^{2} + (0,1)^{2}$ & $T_1$ & Yes & Yes \\
					& $L_{134}(\#0;\#1)$ & $p \neq 2$ & $(1,T(4;0))^2 + (0,T(4;0))^3 + (2,0) + (0,2) + (1,0)^2$ & $A_1$ & Yes & Yes \\
					& $L_{23}(\#1)$ & $p \neq 2$ & $(2,0) + (0,2) + (1,1)^{5} + (1,0)^{4} + (0,1)^{4} + 0^{10}$ & $B_{2}$ & Yes & Yes \\
					& & $p = 2$ & $(0 | ((2,0) + (0,2)) | 0^{2}) + (1,1)^{5} + (1,0)^{4} + (0,1)^{4} + 0^{9}$ & $B_{2}$ & No & No \\
					& $L_{23}(\#2)$ & $p = 2$ & $( (0^{2} + (2,2)) | ( (2,0) + (0,2) ) | 0) + ( ( (0,2) + (2,0) ) | 0)^{4} + (1,1)^{4} + 0^{5}$ & $B_{2}$ & Yes & No \\
					& $L_{123}(\#5)$ & $p \neq 2$, $r \neq 0$ & $(1 \otimes 1^{[r]},2) + (1 \otimes 1^{[r]},0)^{2} + (1^{[r]},1)^{2} + (1,1)^{2} + (2,0) + (2^{[r]},0) + (0,2)^{3} + 0$ & $T_{1}$ & Yes & Yes \\
					& & $p \neq 2$, $r = 0$ & $(2,2) + (2,0)^{4} + (0,2)^{4} + (1,1)^{4} + 0^3$ & $\tilde{A}_{1}$ & Yes & Yes \\
					& & $p = 2$, $r \neq 0$ & $(1 \otimes 1^{[r]},W(2)) + (1 \otimes 1^{[r]},0)^{2} + (1^{[r]},1)^{2} + (1,1)^{2} + (0^{2}|(2,0)+(2^{[r]},0)+(0,2)|0^{2}) + (0,W(2))^{2}$ & $\tilde{A}_{1}$ & No & No \\
					& $L_{123}(\#6)$ & $p \neq 2$ & $(1^{[r]} \otimes 2^{[s]},1) + (1^{[r]},1)^{2} + (1^{[s]},1)^{2} + ({1^{[r]} \otimes 1^{[s]}},0)^{2} + (2^{[r]},0) + (0,2) + (2^{[s]},0)^{3} + 0$ & $T_{1}$ & Yes & Yes \\
					& & $p = 2$ & $(1^{[r]} \otimes W(2)^{[s]},1) + (1^{[r]},1)^{2} + (1^{[s]},1)^{2} + ({1^{[r]} \otimes 1^{[s]}},0)^{2} + (0^{2}|((2^{[r]},0)+(2^{[s]},0)+(0,2))|0^{2}) + (W(2)^{[s]},0)^{2}$ & $\tilde{A}_{1}$ & No & No \\
					& $L_{123}(\#11)$ & $p = 2$ & $(((2^{[r]},2)+(2 \otimes 2^{[r]},0)+(2,2) + 0^{3}) | ((2,0) + (2^{[r]},0) + (0,2)) | 0) + (((2,0)+(2^{[r]},0)+(0,2)) | 0)^{2} + (1 \otimes 1^{[r]},1)^{2}$ & $\tilde{A}_{1}$ & Yes & No \\
					& $L_{234}(\#2)$ & $p \ge 5$ & $(4,2) + (4,1)^2 + (2,0) + (0,3)^2 + (0,2) + 0^{3}$ & $\bar{A}_{1}$ & Yes & Yes \\
					& & $p = 3$ & $(T(4),2) + (2,0) + ( (0,1) | ( (4,1) + (0,3) ) | (0,1))^2 + 0^{3}$ & $\bar{A}_{1}$ & Yes & No \\
					& $L_{234}(\#7)$ & $p \ge 5$ & $(2,0) + (0,T(6;2)) + (0,3)^{2} + (1,3) + (1,4)^{2} + 0^{3}$ & $\bar{A}_{1}$ & Yes & $p \neq 5$ \\
					& $L_{234}(\#9)$ & $p = 2$, $r \neq s, s+1$ & $( 0^{2} + (2^{[s]},2) | ((2^{[r]},0) + (2^{[s]},0) + (0,2)) | 0^{2}) + (1^{[s]},1)^{2} + (1^{[r]} \otimes 1^{[s]},1) + (((1^{[r]} \otimes 2^{[s]},0) + (1^{[r]},2)) | (1^{[r]},0))^{2} + 0^{2}$ & $\bar{A}_{1}$ & No & No \\
					& & $p = 2$, $r=s$ & $( 0 + (2,2) | ((2,0) + (2,0)) | 0) + (0,T(2)^{[1]}) + (1,1)^{2} + (T(2),1) + (3,0)^{2} + (1,W(2))^{2} + 0^{2}$ & $\bar{A}_{1}$ & No & No \\
					& & $p = 2$, $r=s+1$ & $( 0^{2} + (2,2) | ((4,0) + (2,0) + (0,2)) | 0^{2}) + (1,1)^{2} + (3,1) + (((T(2)^{[1]},0) + (2,2)) | (2,0))^{2} + 0^{2}$ & $\bar{A}_{1}$ & No & No \\
					& $L_{234}(\#10)$ & $p = 2$ & $((0^{2} + (0,2 \otimes 2^{[r]})) | ((2,0) + (0,2) + (0,2^{[r]})) | 0^{2}) + {(0,1 \otimes 1^{[r]})^{2}} + {(1,1 \otimes 1^{[r]})} +  (((1,2) + (1,2^{[r]}))|(1,0))^{2} + 0^{2}$ & $\bar{A}_{1}$ & No & No \\
					& $L_{234}(\#13)$ & $p \neq 2$ & $(2,0) + (2^{[r]},0) + (0,2) + (1 \otimes 1^{[r]},1)^{2} + (1,1) + (1^{[r]},1) + (1 \otimes 1^{[r]},0) + (1,0)^{2} + (1^{[r]},0)^{2} + (0,1)^{2} + 0^{3}$ & $\bar{A}_{1}$ & Yes & Yes \\
					& & $p = 2$ & $( 0 | ( (2,0) + (2^{[r]},0) + (0,2) ) | 0^{3}) + (1 \otimes 1^{[r]},1)^{2} + (1,1) + (1^{[r]},1) + (1 \otimes 1^{[r]},0) + (1,0)^{2} + (1^{[r]},0)^{2} + (0,1)^{2} + 0^{2}$ & $\bar{A}_{1}$ & No & No \\ \hline
$A_{1}$				& $L_{1}(\#0)$ & & $T(2;0) + 1^{14} + 0^{20}$ & $C_3$ & Yes & $p \neq 2$ \\
					& $L_{3}(\#0)$ & $p \neq 2$ & $2^7 + 1^8 + 0^{15}$ & $\bar{A}_3$ & Yes & Yes \\
					& 			   & $p = 2$ & $T(2) + W(2)^6 + 1^8 + 0^{14}$ & $B_3$ & Yes & No \\
					& $L_{12}(\#1)$ & $p \neq 2$ & $T(4;0) + 2^{13} + 0^7$ & $\tilde{A}_2$ & Yes & Yes \\
					& $L_{34}(\#1)$ & $p \neq 2$ & $T(4;0)^7 + 2 + 0^7$ & $G_2$ & Yes & Yes \\
					& $L_{13}(\#0^{[r]};\#0^{[s]})$ & $p \neq 2$, $r \neq s$ & $2^{[r]} + (2^{[s]})^{3} + (1^{[r]})^{4} + (1^{[s]})^{4} + (1^{[r]} \otimes 1^{[s]})^{2} + (1^{[r]} \otimes 2^{[s]})^{2} + 0^{4}$ & $\bar{A}_1 T_1$ & Yes & Yes \\
					& & $p \neq 2$, $r=s$ & $2^{6} + 1^{8} + T(3;1)^{2} + 0^{6}$ & $\bar{A}_1 A_1$ & Yes & Yes \\
					& & $p = 2$, $r \neq s, s+1$ & $T(2)^{[r]} + T(2)^{[s]} + (1^{[r]} \otimes W(2)^{[s]})^{2} + (1^{[r]})^{4} + (1^{[s]})^{4} + (1^{[r]} \otimes 1^{[s]})^{2} + (W(2)^{[s]})^{2} + 0^{2}$ & $\bar{A}_1 \tilde{A}_1$ & Yes & No \\
					& & $p = 2$, $r = s+1$ & $T(2) + T(2)^{[1]} + (0|4|0|2)^{2} + 1^{4} + 2^{4} + 3^{2} + W(2)^{2} + 0^{2}$ & $\bar{A}_1 \tilde{A}_1$ & Yes & No \\
					& & $p = 2$, $r=s$ & $3^{2} + T(2)^{4} + W(2)^{2} + 1^{10} + 0^{2}$ & $\tilde{A}_{2}$ & Yes & No \\
					& $L_{23}(\#3)$ & $p \ge 5$ & $T(6;2) + 4^{4} + 3^{4} + 0^{6}$ & $\bar{A}_{1}^2$ & Yes & $p \neq 5$ \\
					& & $p = 2$ & $( (0^{2} + 6) | ( 2 + 4 ) | 0) + ((2 + 4) | 0)^{4} + 3^{4} + 0^{5}$ & $B_{2}$ & Yes & No \\
					& $L_{23}(\#4)$ & $p \neq 2$ & $2 + 2^{[r]} + (1 \otimes 1^{[r]})^{5} + 1^{4} + (1^{[r]})^{4} + 0^{10}$ & $B_{2}$ & Yes & Yes \\
					& & $p = 2$ & $(0 | (2^{[r]} + 2) | 0^{2}) + (1 \otimes 1^{[r]})^{5} + 1^{4} + (1^{[r]})^{4} + 0^{9}$ & $B_{2}$ & No & No \\
					& $L_{23}(\#5)$ & $p = 2$ & $( (0^{2} + (2,2^{[r]})) | ( 2 + 2^{[r]} ) | 0) + ( ( 2 + 2^{[r]} ) | 0)^{4} + ({1 \otimes 1^{[r]}})^{4} + 0^{5}$ & $B_{2}$ & Yes & No \\
					& $L_{134}(\#0^{[r]};\#1^{[s]})$ & $p \neq 2$, $r \neq s, s+1$ & $(1^{[r]} \otimes T(4;0)^{[s]})^{2} + (T(4;0)^{[s]})^3 + 2^{[r]} + 2^{[s]} + (1^{[r]})^2$ & $A_1$ & Yes & Yes \\
					& & $p \neq 2$, $r=s$ & $T(5;3;1)^{2} + 1^{2} + 2^{2} + T(4;0)^3$ & $A_1$ & Yes & Yes \\
					& & $p = 3$, $r = s+1$ & $(1^{[1]} | ( (1 \otimes 2^{[1]}) + 1 ) | 1^{[1]})^{2} + T(4;0)^3 + 2 + 2^{[1]} + (1^{[1]})^2$ & $A_1$ & Yes & Yes \\
					& & $p \geq 5$, $r =s+1$ & $(4 \otimes 1^{[1]})^{2} + (1^{[1]})^{4} + 2 + 2^{[1]} + 4^{3} + 0^{3}$ & $A_1$ & Yes & Yes \\
					& $L_{124}(\#1^{[r]};\#0^{[s]})$ & $p \neq 2$, $r \neq s$ & $T(4;0)^{[s]} + (2^{[s]})^{3} + 2^{[r]} + (1^{[r]} \otimes 2^{[s]})^{2} + ({2^{[r]} \otimes 2^{[s]}})^{2} + (1^{[r]})^{2}$ & $T_{1}$ & Yes & Yes \\
					& & $p \neq 2$, $r=s$ & $T(4;0)^{3} + 2^{6} + T(3;1)^{2} + 1^{2}$ & $A_{1}$ & Yes & Yes \\
					& $L_{123}(\#7)$ &$p \neq 2$, $r,s,t$~distinct & $(1^{[r]} \otimes 1^{[s]} \otimes 2^{[t]}) + 2^{[r]} + 2^{[s]} + (2^{[t]})^{3} + (1^{[r]} \otimes 1^{[s]})^{2} + (1^{[s]} \otimes 1^{[t]})^{2} + (1^{[r]} \otimes 1^{[t]})^{2} + 0$ & $T_{1}$ & Yes & Yes \\
					& & $p \neq 2$, $r = s < t$ & $(2 \otimes 2^{[t]}) + 2^{4} + (2^{[t]})^{4} + (1 \otimes 1^{[t]})^{4} + 0^{3}$ & $\tilde{A}_{1}$ & Yes & Yes \\
					& & $p \neq 2$, $r = t < s$ & $(T(3;1) \otimes 1^{[s]}) + 2^{6} + 2^{[s]} + (1 \otimes 1^{[s]})^{4} + 0^{3}$ & $\tilde{A}_{1}$ & Yes & Yes \\
					& & $p \neq 2$, $r < s = t$ & $2 + (2^{[s]})^{6} + (1 \otimes T(3;1)^{[s]}) + (1 \otimes 1^{[s]})^{4} + 0^{3}$ & $\tilde{A}_{1}$ & Yes & Yes \\
					& & $p = 2$, $r,s,t$~distinct & $(1^{[r]} \otimes 1^{[s]} \otimes W(2)^{[t]}) + (1^{[r]} \otimes 1^{[s]})^{2} + (1^{[s]} \otimes 1^{[t]})^{2} + (1^{[r]} \otimes 1^{[t]})^{2} +~(W(2)^{[t]})^{2} +(0^{2}|(2^{[r]} + 2^{[s]} + 2^{[t]})|0^{2})$ & $\tilde{A}_{1}$ & No & No \\
					& & $p = 2$, $r < s = t$ & $(1 \otimes 3^{[s]}) + (1 \otimes 1^{[s]})^{5} + (T(2)^{[s]})^{3} + {(0 | (2 + 2^{[s]}) | 0)} + (W(2)^{[s]})^{2}$ & $\tilde{A}_{1}$ & No & No \\
					& & $p = 2$, $r = t < s$ & $(3 \otimes 1^{[s]}) + (1 \otimes 1^{[s]})^{5}+ T(2)^{3} + (0 | (2 + 2^{[s]})| 0) + W(2)^{2}$ & $\tilde{A}_{1}$ & No & No \\
					& $L_{123}(\#9)$ & $p \geq 7$ & $T(10;2) + 6^{5} + 0^{3}$ & $A_{1}$ & Yes & $p \neq 7$ \\
					& $L_{123}(\#12)$ & $p = 2$ & $(( (2^{[r]} \otimes 2^{[s]}) + (2 \otimes 2^{[r]})+(2 \otimes 2^{[s]}) + 0^{3}) | (2 + 2^{[r]} + 2^{[s]}) | 0) + ((2+2^{[r]}+2^{[s]}) | 0)^{2} + (1 \otimes 1^{[r]} \otimes 1^{[s]})^{2}$ & $\tilde{A}_{1}$ & Yes & No \\
					& $L_{234}(\#5)$ & $p \geq 7$ & $T(10;2) + 6 + T(9;3)^{2} + 0^{3}$ & $\bar{A}_{1}$ & Yes & $p \neq 7$ \\
					& $L_{234}(\#11)$
						& $p = 2$, $r\neq s,t,s+1,t+1$ & $((0^{2} + (2^{[s]} \otimes 2^{[t]})) | (2^{[r]} + 2^{[s]} + 2^{[t]}) | 0^{2}) + ({1^{[s]} \otimes 1^{[t]}})^{2} + ({1^{[r]} \otimes 1^{[s]} \otimes 1^{[t]}}) +  (((1^{[r]} \otimes 2^{[s]}) + (1^{[r]} \otimes 2^{[t]}))| 1^{[r]} )^{2} + 0^{2}$ & $\bar{A}_{1}$ & No & No \\
						& & $p = 2$, $r = s$, $t \ge 2$ & $((0 + (2 \otimes 2^{[t]})) | (2 + 2^{[t]}) | 0) + T(2) + (1 \otimes 1^{[t]})^{2} + (T(2) \otimes 1^{[t]}) + 3^{2} + ( (1 \otimes 2^{[t]}) | 1 )^{2} + 0^{2}$ & $\bar{A}_{1}$ & No & No \\
						& & $p = 2$, $r = s$, $t = 1$ & $((0 + 6) | (2 + 4) | 0) + T(2) + 3^{4} + T(4) + W(5)^{2} + 0^{2}$ & $\bar{A}_{1}$ & No & No \\
						& & $p = 2$, $r=s+1$, $t \geq 2$ & $((0^{2} + (2 \otimes 2^{[t]})) | (2 + 4 + 2^{[t]}) | 0^{2}) + (1 \otimes 1^{[t]})^{2} + ({3 \otimes 1^{[t]}}) +  (0|4|(0 + (2 \otimes 2^{[t]}))| 2 )^{2} + 0^{2}$ & $\bar{A}_{1}$ & No & No \\
						& & $p = 2$, $r = t+1$ & $((0^{2} + (2 \otimes 4^{[t]})) | (2 + 2^{[t]} + 4^{[t]}) | 0^{2}) + (1 \otimes 1^{[t]})^{2} + (1 \otimes 1^{[t]} \otimes 2^{[t]}) +  (0|4^{[t]}|((2 \otimes 2^{[t]}) + 0)| 2^{[t]} )^{2} + 0^{2}$ & $\bar{A}_{1}$ & No & No \\
						& & $p = 2$, $r=t=1$ & $((0 + (2 \otimes 4)) | (2 + 4) | 0) + T(2)^{[1]} + 3^{2} + ({T(2)^{[1]} \otimes 1}) + (2 \otimes 4)^2 + (( 0 | 4 | 0| 2 )^{2} + 0^{2}$ & $\bar{A}_{1}$ & No & No \\
						& & $p = 2$, $r=t > 1$ & $((0 + (2 \otimes 2^{[r]})) | (2 + 2^{[r]}) | 0) + T(2)^{[r]} + ({1 \otimes 1^{[r]}})^{2} + (1 \otimes T(2)^{[r]}) + (1^{[r]} \otimes 2^{[r]})^2 + ((2 \otimes 1^{[r]}) | 1^{[r]})^{2} + 0^{2}$ & $\bar{A}_{1}$ & No & No \\
					& $L_{234}(\#12)$ & $p \ge 5$, $r \neq s$ & $2^{[r]} + T(6;2)^{[s]} + (3^{[r]})^{2} + ({1^{[r]} \otimes 3^{[s]}}) + ({1^{[r]} \otimes 4^{[s]}})^{2} + 0^{3}$ & $\bar{A}_{1}$ & Yes & $p\neq 5$ or $r < s$ \\   
					& & $p \ge 5$, $r = s$ & $2^{2} + T(6;2) + 4 + 3^{2} + T(5;3)^{2} + 0^{3}$ & $\bar{A}_{1}$ & Yes & Yes \\
					& $L_{234}(\#14)$ & $p \neq 2$ & $2+2^{[r]}+2^{[s]} + (1 \otimes 1^{[r]} \otimes 1^{[s]})^2 + (1 \otimes 1^{[r]}) + ({1 \otimes 1^{[s]}}) + (1^{[r]} \otimes 1^{[s]}) + 1^{2} + (1^{[r]})^{2} + (1^{[s]})^{2} + 0^{3}$ & $\bar{A}_{1}$ & Yes & Yes \\
					& & $p = 2$ & $(0 | (2 + 2^{[r]} + 2^{[s]}) | 0^{3}) + (1 \otimes 1^{[r]} \otimes 1^{[s]})^{2} + ({1 \otimes 1^{[r]}}) + (1 \otimes 1^{[s]}) + (1^{[r]} \otimes1^{[s]}) + 1^{2} + (1^{[r]})^{2} + (1^{[s]})^{2} + 0^{2}$ & $\bar{A}_{1}$ & No & No \\
					& $L_{234}(\#15)$ & $p \ge 5$ & $(4^{[r]} \otimes 2^{[s]}) + (4^{[r]} \otimes 1^{[s]})^2 + 2^{[r]} + (3^{[s]})^2 + 2^{[s]} + 0^{3}$ & $\bar{A}_{1}$ & Yes & Yes \\
					& & $p = 3$, $s \neq r+1$ & $(T(4)^{[r]} \otimes 2^{[s]}) + 2^{[r]} + ( 1^{[s]} | ( (4^{[r]} \otimes 1^{[s]}) + 3^{[s]} ) | 1^{[s]})^2 + 0^{3}$ & $\bar{A}_{1}$ & Yes  & $r < s$  \\
					& & $p = 3$, $s = r+1$ & $T(10) + 2 + ( 3 | (9 + 7 +1) | 3)^2 + 0^{3}$ & $\bar{A}_{1}$ & Yes  & Yes  \\
					\hline
\end{longtable}
}

\bibliographystyle{amsplain}
\bibliography{biblio}

\providecommand{\bysame}{\leavevmode\hbox to3em{\hrulefill}\thinspace}
\providecommand{\MR}{\relax\ifhmode\unskip\space\fi MR }
\providecommand{\MRhref}[2]{%
  \href{http://www.ams.org/mathscinet-getitem?mr=#1}{#2}
}
\providecommand{\href}[2]{#2}
\begin{thebibliography}{10}

\bibitem{MR2707891}
Bonnie Amende, \emph{G-irreducible subgroups of type {$A_1$}}, ProQuest LLC,
  Ann Arbor, MI, 2005, Thesis (Ph.D.)--University of Oregon. \MR{2707891}

\bibitem{MR898346}
Michael Aschbacher, \emph{Chevalley groups of type {$G_2$} as the group of a
  trilinear form}, J. Algebra \textbf{109} (1987), no.~1, 193--259. \MR{898346
  (88g:20089)}

\bibitem{MR1047327}
H.~Azad, M.~Barry, and G.~Seitz, \emph{On the structure of parabolic
  subgroups}, Comm. Algebra \textbf{18} (1990), no.~2, 551--562. \MR{1047327
  (91d:20048)}

\bibitem{MR2178661}
Michael Bate, Benjamin Martin, and Gerhard R{\"o}hrle, \emph{A geometric
  approach to complete reducibility}, Invent. Math. \textbf{161} (2005), no.~1,
  177--218. \MR{2178661 (2007k:20101)}

\bibitem{MR0294349}
A.~Borel and J.~Tits, \emph{\'{E}l\'ements unipotents et sous-groupes
  paraboliques de groupes r\'eductifs. {I}}, Invent. Math. \textbf{12} (1971),
  95--104. \MR{0294349 (45 \#3419)}

\bibitem{MR0240238}
N.~Bourbaki, \emph{\'{E}l\'ements de math\'ematique. {F}asc. {XXXIV}. {G}roupes
  et alg\`ebres de {L}ie. {C}hapitre {IV}: {G}roupes de {C}oxeter et syst\`emes
  de {T}its. {C}hapitre {V}: {G}roupes engendr\'es par des r\'eflexions.
  {C}hapitre {VI}: syst\`emes de racines}, Actualit\'es Scientifiques et
  Industrielles, No. 1337, Hermann, Paris, 1968. \MR{0240238 (39 \#1590)}

\bibitem{MR0407163}
Roger~W. Carter, \emph{Simple groups of {L}ie type}, John Wiley \& Sons,
  London-New York-Sydney, 1972, Pure and Applied Mathematics, Vol. 28.
  \MR{0407163 (53 \#10946)}

\bibitem{MR1132853}
Arjeh~M. Cohen, Martin~W. Liebeck, Jan Saxl, and Gary~M. Seitz, \emph{The local
  maximal subgroups of exceptional groups of {L}ie type, finite and algebraic},
  Proc. London Math. Soc. (3) \textbf{64} (1992), no.~1, 21--48. \MR{1132853
  (92m:20012)}

\bibitem{Craven2017}
David~A. Craven, \emph{Alternating subgroups of exceptional groups of {L}ie
  type}, Proceedings of the London Mathematical Society \textbf{115} (2017),
  no.~3, 449--501.

\bibitem{DonkinTilting}
Stephen Donkin, \emph{On tilting modules for algebraic groups}, Mathematische
  Zeitschrift \textbf{212} (1993), no.~1, 39--60 (English).

\bibitem{Dot1}
Stephen Doty, \emph{Weyl modules}, Gap Package, ver. 1.1.

\bibitem{MR1490581}
Daniel Gorenstein, Richard Lyons, and Ronald Solomon, \emph{The classification
  of the finite simple groups. {N}umber 3. {P}art {I}. {C}hapter {A}},
  Mathematical Surveys and Monographs, vol.~40, American Mathematical Society,
  Providence, RI, 1998, Almost simple $K$-groups. \MR{1490581 (98j:20011)}

\bibitem{MR1466951}
R.~Lawther and D.~M. Testerman, \emph{{$A_1$} subgroups of exceptional
  algebraic groups}, Mem. Amer. Math. Soc. \textbf{141} (1999), no.~674,
  viii+131. \MR{1466951}

\bibitem{MR1329942}
Martin~W. Liebeck and Gary~M. Seitz, \emph{Reductive subgroups of exceptional
  algebraic groups}, Mem. Amer. Math. Soc. \textbf{121} (1996), no.~580,
  vi+111. \MR{1329942 (96i:20059)}

\bibitem{MR1650328}
\bysame, \emph{On the subgroup structure of classical groups}, Invent. Math.
  \textbf{134} (1998), no.~2, 427--453. \MR{1650328 (99h:20074)}

\bibitem{MR2044850}
\bysame, \emph{The maximal subgroups of positive dimension in exceptional
  algebraic groups}, Mem. Amer. Math. Soc. \textbf{169} (2004), no.~802,
  vi+227. \MR{2044850 (2005b:20082)}

\bibitem{MR2043006}
Martin~W. Liebeck and Donna~M. Testerman, \emph{Irreducible subgroups of
  algebraic groups}, Q. J. Math. \textbf{55} (2004), no.~1, 47--55. \MR{2043006
  (2005b:20087)}

\bibitem{Litterick2018}
Alastair~J. Litterick, \emph{On non-generic finite subgroups of exceptional
  algebraic groups}, Mem. Amer. Math. Soc. \textbf{253} (2018), no.~1207,
  vi+156.

\bibitem{LitTho}
Alastair~J. Litterick and Adam~R. Thomas, \emph{Complete reducibility in good
  characteristic}, Trans. Amer. Math. Soc. \textbf{370} (2018), no.~8,
  5279--5340.

\bibitem{LitTho2}
\bysame, \emph{Complete reducibility in bad characteristic {I}}, in
  preparation.

\bibitem{MR1901354}
Frank L{\"u}beck, \emph{Small degree representations of finite {C}hevalley
  groups in defining characteristic}, LMS J. Comput. Math. \textbf{4} (2001),
  135--169 (electronic). \MR{1901354 (2003e:20013)}

\bibitem{MR2850737}
Gunter Malle and Donna Testerman, \emph{Linear algebraic groups and finite
  groups of {L}ie type}, Cambridge Studies in Advanced Mathematics, vol. 133,
  Cambridge University Press, Cambridge, 2011. \MR{2850737 (2012i:20058)}

\bibitem{Richardson1977}
R.~W. Richardson, \emph{Affine coset spaces of reductive algebraic groups}, The
  Bulletin of the London Mathematical Society \textbf{9} (1977), no.~1, 38--41.
  \MR{0437549}

\bibitem{MR888704}
Gary~M. Seitz, \emph{The maximal subgroups of classical algebraic groups}, Mem.
  Amer. Math. Soc. \textbf{67} (1987), no.~365, iv+286. \MR{888704 (88g:20092)}

\bibitem{Ser3}
Jean-Pierre Serre, \emph{Morsund lectures, {U}niversity of {O}regon}, 1998.

\bibitem{steinberglectureson}
Robert Steinberg, \emph{Lectures on {C}hevalley groups}, Notes prepared by J.
  Faulkner and R. Wilson, Yale University, New Haven, Conn., 1968.

\bibitem{MR2604850}
David~I. Stewart, \emph{The reductive subgroups of {$G_2$}}, J. Group Theory
  \textbf{13} (2010), no.~1, 117--130. \MR{2604850 (2011c:20099)}

\bibitem{MR3075783}
\bysame, \emph{The reductive subgroups of {$F_4$}}, Mem. Amer. Math. Soc.
  \textbf{223} (2013), no.~1049, vi+88. \MR{3075783}

\bibitem{Tho1}
Adam~R. Thomas, \emph{Simple irreducible subgroups of exceptional algebraic
  groups}, J. Algebra \textbf{423} (2015), 190--238.

\bibitem{Tho2}
\bysame, \emph{Irreducible {$A_1$} subgroups of exceptional algebraic groups},
  J. Algebra \textbf{447} (2016), 240--296.

\bibitem{Tho3}
\bysame, \emph{The irreducible subgroups of exceptional algebraic groups}, to
  appear in Mem. Amer. Math. Soc.

\end{thebibliography}

\end{document}